\documentclass{cocv}

\usepackage{amssymb} %in order to use, for example, $\varkappa$.

\def\ess{{\rm ess}} %for the esssup
\def\C{{\rm C}} %for the space of continuous functions
\def\AC{{\rm AC}} %for the space of absolutely continuous functions
\def\Lip{{\rm Lip}} %for the space of Lipschitz continuous functions
\def\dist{{\rm dist}} %for the distance

\begin{document}

\title{Minimax Solutions of Hamilton--Jacobi Equations with Fractional Coinvariant Derivatives}
\thanks{This work was supported by RSF (project no.~19-71-00073).}

\author{Mikhail Gomoyunov}
\address{N.N.~Krasovskii Institute of Mathematics and Mechanics of the Ural Branch of the Russian Academy of Sciences;}
\secondaddress{Ural Federal University; \email{m.i.gomoyunov@gmail.com}}

\date{The dates will be set by the publisher.}

\begin{abstract}
    We consider a Cauchy problem for a Hamilton--Jacobi equation with coinvariant derivatives of an order $\alpha \in (0, 1)$.
    Such problems arise naturally in optimal control problems for dynamical systems which evolution is described by ordinary differential equations with the Caputo fractional derivatives of the order $\alpha$.
    We propose a notion of a generalized in the minimax sense solution of the considered problem.
    We prove that a minimax solution exists, is unique, and is consistent with a classical solution of this problem.
    In particular, we give a special attention to the proof of a comparison principle, which requires construction of a suitable Lyapunov--Krasovskii functional.
\end{abstract}

%\begin{resume}
%
%\end{resume}



\maketitle

\section*{Introduction}
    Nowadays, the theory of differential equations with fractional-order derivatives (see, e.g., \cite{Miller_Ross_1993,Samko_Kilbas_Marichev_1993,Podlubny_1999,Kilbas_Srivastava_Trujillo_2006,Diethelm_2010}) is an actively developing branch of mathematics, which attracts the interest of many researchers.
    In particular, attention is paid to optimal control problems for dynamical systems which evolution is described by ordinary differential equations with the Caputo fractional derivatives.
    Such problems appear in various fields of knowledge including, e.g., chemistry \cite{Flores-Tlacuahuac_Biegler_2014}, biology \cite{Toledo-Hernandez_et_all_2014}, electrical engineering \cite{Kaczorek_2016}, and medicine \cite{Kheiri_Jafari_2018}.
    Main directions of research here are related to necessary optimality conditions (see, e.g., \cite{Bergounioux_Bourdin_2020,Lin_Yong_2020} and the references therein) and numerical methods for constructing optimal controls (see, e.g., \cite{Zeid_Effati_Kamyad_2018,Li_Wang_Rehbock_2019,Salati_Shamsi_Torres_2019} and the references therein).
    In addition, note that several problems for linear systems are considered and studied in detail in, e.g., \cite{Kamocki_Majewski_2014,Kubyshkin_Postnov_2014_Eng,Idczak_Walczak_2016,Matychyn_Onyshchenko_2018,Bandaliyev_et_all_2020}.
    The reader is also referred to \cite{Butkovskii_Postnov_Postnova_2013_Eng} for an overview of works on various control problems for fractional-order systems.

    In \cite{Gomoyunov_2020_SIAM}, the dynamic programming principle was extended to a Bolza-type optimal control problem for a dynamical system described by a fractional differential equation with the Caputo derivative of an order $\alpha \in (0, 1)$.
    In particular, it was shown that the value of this problem should be introduced as a functional in a suitable space of paths.
    Further, the problem was associated with a Hamilton--Jacobi equation with coinvariant ($ci$-) derivatives of the order $\alpha$.
    Note that these derivatives can be considered as a suitable extension of the notion of $ci$-derivatives (of the first order) proposed and developed in, e.g., \cite{Kim_1999,Lukoyanov_2011_Eng}.
    It was proved that if the value functional is smooth enough (namely, if it is $ci$-smooth of the order $\alpha$), then it satisfies the Hamilton--Jacobi equation and the natural boundary condition, and, therefore, the value functional can be treated as a solution of this Cauchy problem in the classical sense.
    However, as a rule, the value functional does not possess the required smoothness properties, which leads to the need to introduce and study generalized solutions of the obtained Cauchy problem.

    In the paper, we consider a Cauchy problem for a Hamilton--Jacobi equation with $ci$-derivatives of an order $\alpha \in (0, 1)$ and propose a notion of a minimax solution of this problem.
    The technique of minimax solutions originates in the positional differential games theory (see, e.g., \cite{Krasovskii_Subbotin_1988,Krasovskii_Krasovskii_1995}) and can be seen as the development of the unification constructions of differential games \cite{Krasovskii_1976}.
    Minimax solutions of Hamilton--Jacobi equations with first-order partial derivatives were proposed and comprehensively studied in \cite{Subbotin_1995} (see also \cite{Subbotin_1996}).
    Further, this technique was extended to Hamilton--Jacobi equations with first-order $ci$-derivatives, which arise in optimization problems for dynamical systems described by functional differential equations of a retarded type \cite{Lukoyanov_2011_Eng} (see also \cite{Lukoyanov_2003_1,Lukoyanov_2003_2,Lukoyanov_2004_PMM_Eng,Lukoyanov_2006_IMM_Eng}, and \cite{Bayraktar_Keller_2018} for an infinite dimensional case) and of a neutral type \cite{Lukoyanov_Gomoyunov_Plaksin_2017_Doklady,Lukoyanov_Plaksin_2019_MIAN_Eng,Plaksin_2019_DE_Eng}.
    Note that the minimax approach was also applied to investigate generalized solutions of systems of equations arising in mean field games \cite{Averboukh_2015}.

    Following the general methodology, we define a minimax solution of the considered Cauchy problem in terms of a pair of non-local stability properties of this solution with respect to so-called characteristic differential inclusions, which in this case become fractional differential inclusions with the Caputo derivatives of the order $\alpha$.
    We prove that a minimax solution exists, is unique, and is consistent with a classical solution of the problem.
    In particular, we establish a comparison principle.
    In general, the proofs of these results are carried out by the schemes of the proofs of the corresponding statements for Hamilton--Jacobi equations with partial derivatives \cite{Subbotin_1995} and with first-order $ci$-derivatives \cite{Lukoyanov_2011_Eng} (see also \cite{Lukoyanov_2003_1} and  \cite{Bayraktar_Keller_2018}).
    They are based on properties \cite{Gomoyunov_2020_DE} of the sets of solutions of the characteristic differential inclusions.
    However, in order to prove the comparison principle, it is required to construct a suitable Lyapunov--Krasovskii functional with a number of prescribed properties (in this connection, see, e.g., \cite[Sect.~15]{Lukoyanov_2011_Eng} and also \cite[Sect.~5]{Lukoyanov_2004_PMM_Eng}).
    Due to features of fractional-order integrals and derivatives, this construction substantially differs from the previous studies and can be considered as the main contribution of the paper.

    The paper is organized as follows.
    In Sect.~\ref{section_Preliminaries}, we recall definitions of Riemann--Liouville integrals and Caputo derivatives of a fractional order, describe some of their properties, and introduce special functional spaces.
    Auxiliary facts from the theory of differential inclusions with the Caputo fractional derivatives are presented in Sect.~\ref{section_Differential_Inclusions}.
    In Sect.~\ref{section_HJ}, a Cauchy problem for a Hamilton--Jacobi equation with $ci$-derivatives of a fractional order is considered, and a definition of a minimax solution of this problem is given.
    Consistency of minimax and classical solutions of the problem is discussed in Sect.~\ref{section_Consistency}.
    A comparison principle is established in Sect.~\ref{section_Comparison}.
    Existence and uniqueness of a minimax solution are proved in Sect.~\ref{section_Existence_Uniqueness}.
    Concluding remarks are given in Sect.~\ref{section_Conclusion}.

\section{Preliminaries}
\label{section_Preliminaries}

    Fix $n \in \xN$ and $T > 0$.
    By $\|\cdot\|$ and $\langle \cdot, \cdot \rangle$, we denote the Euclidian norm and the inner product in $\xR^n$.

    For a given $t \in [0, T]$, let $\xLinfty([0, t], \xR^n)$ be the set of (Lebesgue) measurable and essentially bounded functions from $[0, t]$ to $\xR^n$.
    For a function $\psi(\cdot) \in \xLinfty([0, t], \xR^n)$, the (left-sided) Riemann--Liouville fractional integral of an order $\alpha > 0$ is defined by
    \begin{equation} \label{I}
        (I^\alpha \psi)(\tau)
        = \frac{1}{\Gamma(\alpha)} \int_{0}^{\tau} \frac{\psi(\xi)}{(\tau - \xi)^{1 - \alpha}} \xdif \xi,
        \quad \tau \in [0, t],
    \end{equation}
    where $\Gamma$ is the gamma function.
    In the case $\alpha = 0$, we formally define $(I^0 \psi)(\tau) = \psi(\tau)$, $\tau \in [0, t]$.

    Note that, for every $\alpha \geq 0$, $\beta \geq 0$, and $\psi(\cdot) \in \xLinfty([0, t], \xR^n)$, the following semigroup property holds (see, e.g., \cite[(2.21)]{Samko_Kilbas_Marichev_1993} and also \cite[Theorem~2.2]{Diethelm_2010}):
    \begin{equation} \label{I_semigroup_property}
        \big(I^\alpha (I^\beta \psi) \big)(\tau)
        = (I^{\alpha + \beta} \psi)(\tau),
        \quad \tau \in [0, t].
    \end{equation}

    Further, according to, e.g., \cite[Theorem~3.6 and Remark~3.3]{Samko_Kilbas_Marichev_1993} (see also \cite[Theorem~2.6]{Diethelm_2010}), for any $\alpha \in (0, 1]$ and $\psi(\cdot) \in \xLinfty([0, t], \xR^n)$, the inequalities below are valid:
    \begin{equation} \label{I_properties}
        \|(I^\alpha \psi)(\tau)\|
        \leq \frac{\tau^\alpha}{\Gamma(\alpha + 1)} \ess \sup_{\xi \in [0, \tau]} \|\psi(\xi)\|,
        \quad \|(I^\alpha \psi)(\tau) - (I^\alpha \psi)(\tau^\prime)\|
        \leq \frac{2 |\tau - \tau^\prime|^\alpha}{\Gamma(\alpha + 1)} \ess \sup_{\xi \in [0, t]} \|\psi(\xi)\|,
    \end{equation}
    where $\tau$, $\tau^\prime \in [0, t]$.
    In particular, we have the inclusion $(I^\alpha \psi)(\cdot) \in \C([0, t], \xR^n)$.
    Here and below, by $\C([0, t], \xR^n)$, we denote the space of continuous functions from $[0, t]$ to $\xR^n$ endowed with the norm
    \begin{equation*}
        \|x(\cdot)\|_{[0, t]}
        = \max_{\tau \in [0, t]} \|x(\tau)\|,
        \quad x(\cdot) \in \C([0, t], \xR^n).
    \end{equation*}

    For $\alpha \in (0, 1]$, let $\AC^\alpha([0, t], \xR^n)$ be the set of functions $x: [0, t] \to \xR^n$ that can be represented in the form
    \begin{equation} \label{AC^alpha}
        x(\tau)
        = x(0) + (I^\alpha \psi)(\tau),
        \quad \tau \in [0, t],
    \end{equation}
    for some function $\psi(\cdot) \in \xLinfty([0, t], \xR^n)$.
    The set $\AC^\alpha([0, t], \xR^n)$ is considered as a subset of $\C([0, t], \xR^n)$.
    Note that, in the case $\alpha = 1$, the set $\AC^1([0, t], \xR^n)$ coincides with the set $\Lip([0, t], \xR^n)$ of Lipschitz continuous functions from $[0, t]$ to $\xR^n$.

    Let $\alpha \in (0, 1]$ and $x(\cdot) \in \AC^\alpha([0, t], \xR^n)$.
    It follows from \eqref{I_semigroup_property} that, for every $\beta \in [0, 1 - \alpha]$, the inclusion $(I^\beta (x(\cdot) - x(0)))(\cdot) \in \AC^{\alpha + \beta}([0, t], \xR^n)$ holds.
    In particular, we obtain $(I^{1 - \alpha} (x(\cdot) - x(0)))(\cdot) \in \Lip([0, t], \xR^n)$.
    Hence, the (left-sided) Caputo fractional derivative of $x(\cdot)$ of the order $\alpha$, which is defined by
    \begin{equation} \label{Caputo}
        (^C D^\alpha x) (\tau)
        = \xDrv{\big(I^{1 - \alpha} (x(\cdot) - x(0))\big)(\tau)}{\tau},
    \end{equation}
    exists for almost every (a.e.) $\tau \in [0, t]$, and, moreover, the equality $(^C D^\alpha x)(\tau) = \psi(\tau)$ is valid for a.e. $\tau \in [0, t]$, where $\psi(\cdot) \in \xLinfty([0, t], \xR^n)$ is the function from \eqref{AC^alpha}.
    If $\alpha = 1$, then the Caputo derivative $(^C D^1 x)(\tau)$ is the usual first-order derivative $\dot{x}(\tau) = \xDrv{x(\tau)}{\tau}$.

    Now, consider the set $G_n$ of pairs $(t, w(\cdot))$ such that $t \in [0, T]$ and $w(\cdot) \in \C([0, t], \xR^n)$.
    For $x(\cdot) \in \C([0, T], \xR^n)$ and $t \in [0, T]$, let $x_t(\cdot) \in \C([0, t], \xR^n)$ denote the restriction of the function $x(\cdot)$ to the interval $[0, t]$:
    \begin{equation} \label{x_t}
        x_t(\tau)
        = x(\tau),
        \quad \tau \in [0, t].
    \end{equation}
    Then, we have $(t, x_t(\cdot)) \in G_n$.
    In accordance with \cite[\S~1]{Lukoyanov_2011_Eng} (see also, e.g., \cite{Lukoyanov_2003_1}), the set $G_n$ is endowed with the metric
    \begin{equation} \label{dist}
        \dist \big((t, w(\cdot)), (t^\prime, w^\prime(\cdot))\big)
        = \max \big\{ \dist^\ast\big((t, w(\cdot)), (t^\prime, w^\prime(\cdot))\big),
        \dist^\ast\big((t^\prime, w^\prime(\cdot)), (t, w(\cdot))\big) \big\},
    \end{equation}
    where
    \begin{equation*}
        \dist^\ast\big((t, w(\cdot)), (t^\prime, w^\prime(\cdot))\big)
        = \max_{\tau \in [0, t]} \min_{\tau^\prime \in [0, t^\prime]}
        \sqrt{|\tau - \tau^\prime|^2 + \|w(\tau) - w^\prime(\tau^\prime)\|^2},
        \quad (t, w(\cdot)), (t^\prime, w^\prime(\cdot)) \in G_n.
    \end{equation*}

    Let us describe some properties of this metric.
    By \cite[Proposition~8.2]{Gomoyunov_2020_SIAM}, for any $(t, w(\cdot))$, $(t^\prime, w^\prime(\cdot)) \in G_n$, $t^\prime \leq t$, the inequalities
    \begin{equation} \label{dist_properties}
        \dist
        \leq t - t^\prime + \varkappa(t - t^\prime) + \max_{\tau \in [0, t^\prime]} \|w(\tau) - w^\prime(\tau)\|,
        \quad t - t^\prime
        \leq \dist,
        \quad \max_{\tau \in [0, t^\prime]} \|w(\tau) - w^\prime(\tau)\|
        \leq \dist + \varkappa(\dist)
    \end{equation}
    are valid, where $\dist = \dist((t, w(\cdot)), (t^\prime, w^\prime(\cdot)))$ and $\varkappa$ is the modulus of continuity of $w(\cdot)$ on $[0, t]$ given by
    \begin{equation*}
        \varkappa(\delta)
        = \max \big\{ \|w(\tau) - w(\tau^\prime)\|:
        \, \tau, \tau^\prime \in [0, t], \, |\tau - \tau^\prime| \leq \delta \big\},
        \quad \delta
        \geq 0.
    \end{equation*}
    In particular, if sequences $\{x^{(k)}(\cdot)\}_{k \in \xN} \subset \C([0, T], \xR^n)$ and $\{t_k\}_{k \in \xN} \subset [0, T]$ converge to $x^{(0)}(\cdot) \in \C([0, T], \xR^n)$ and $t_0 \in [0, T]$, respectively, then $(t_k, x^{(k)}_{t_k}(\cdot)) \to (t_0, x^{(0)}_{t_0}(\cdot))$ as $k \to \infty$ with respect to the metric $\dist$.

    Moreover, note also that if a sequence $\{(t_k, w_k(\cdot))\}_{k \in \xN} \subset G_n$ is convergent, then the functions $w_k(\cdot)$, $k \in \xN$, are uniformly bounded and equicontinuous (see, e.g., \cite[Assertion~6]{Gomoyunov_2020_DE}).
    Namely, there exists $R > 0$ such that $\|w_k(\cdot)\|_{[0, t_k]} \leq R$, $k \in \xN$, and the function $\varkappa_\ast(\delta) = \sup\{ \varkappa_k(\delta): \, k \in \xN\}$, $\delta \geq 0$, satisfies the relation $\varkappa_\ast(\delta) \to 0$ as $\delta \to + 0$, where $\varkappa_k$ is the modulus of continuity of $w_k(\cdot)$ on $[0, t_k]$.

    Finally, for every $\alpha \in (0, 1]$, we introduce the following two subsets of $G_n$:
    \begin{equation} \label{G^alpha_n}
        G^\alpha_n
        = \big\{ (t, w(\cdot)) \in G_n:
        \, w(\cdot) \in \AC^\alpha([0, t], \xR^n) \big\},
        \quad G_n^{\alpha \circ}
        = \big\{ (t, w(\cdot)) \in G_n^\alpha:
        \, t < T\big\}.
    \end{equation}

\section{Differential Inclusions with Fractional Derivatives}
\label{section_Differential_Inclusions}

    This section provides auxiliary statements concerning properties of the sets of solutions of differential inclusions with the Caputo fractional derivatives of an order $\alpha \in (0, 1)$, which constitute a basis for the proofs of the main results of the paper.

    \subsection{Ordinary Differential Inclusions}
    \label{subsection_ODI}

        Suppose that a set-valued function $[0, T] \times \xR^n \ni (t, x) \mapsto F(t, x) \subset \xR^n \times \xR$ satisfies the following conditions:
        \begin{itemize}
            \item[$(F.1)$]
                For every $t \in [0, T]$ and $x \in \xR^n$, the set $F(t, x)$ is nonempty, convex, and compact in $\xR^n \times \xR$.

            \item[$(F.2)$]
                The set-valued function $F$ is upper semicontinuous (in the Hausdorff sense).
                It means that, for every $(t, x) \in [0, T] \times \xR^n$ and $\varepsilon > 0$, there exists $\delta > 0$ such that, for any $(t^\prime, x^\prime) \in [0, T] \times \xR^n$, the inequality $|t - t^\prime|^2 + \|x - x^\prime\|^2 \leq \delta^2$ implies the inclusion $F(t^\prime, x^\prime) \subset [F(t, x)]^\varepsilon$.
                Here and below, for $\varepsilon > 0$ and $F \subset \xR^n \times \xR$, the symbol $[F]^\varepsilon$ stands for the $\varepsilon$-neighbourhood of $F$, given by
                \begin{equation*}
                    [F]^\varepsilon
                    = \Big\{(f, h) \in \xR^n \times \xR:
                    \, \inf_{(f^\prime, h^\prime) \in F} \big(\|f - f^\prime\|^2 + |h - h^\prime|^2 \big)^{1 / 2}
                    \leq \varepsilon \Big\}.
                \end{equation*}

            \item[$(F.3)$]
                There exists $c_F > 0$ such that, for any $(t, x) \in [0, T] \times \xR^n$,
                \begin{equation*}
                    \sup \big\{\|f\|:
                    \, (f, h) \in F(t, x)\big\}
                    \leq c_F (1 + \|x\|).
                \end{equation*}
        \end{itemize}

        Given $(t_0, w_0(\cdot)) \in G_n^\alpha$ and $z_0 \in \xR$, consider the Cauchy problem for the differential inclusion
        \begin{equation} \label{Differential_Inclusion}
            \big( (^C D^\alpha x)(t), \dot{z}(t) \big)
            \in F(t, x(t)),
            \quad (x(t), z(t)) \in \xR^n \times \xR, \quad t \in [t_0, T],
        \end{equation}
        and the initial condition
        \begin{equation} \label{Differential_Inclusion_Initial_Condition}
            x(t)
            = w_0(t),
            \quad
            z(t)
            = z_0,
            \quad t \in [0, t_0].
        \end{equation}

        Let $XZ^\alpha(t_0, w_0(\cdot), z_0)$ be the set of pairs of functions $(x(\cdot), z(\cdot)) \in \AC^\alpha([0, T], \xR^n) \times \Lip([0, T], \xR)$ satisfying \eqref{Differential_Inclusion_Initial_Condition}.
        Note that it is convenient to identify any such pair $(x(\cdot), z(\cdot))$ with the corresponding function from $[0, T]$ to $\xR^n \times \xR$.
        In this sense, taking into account that $\AC^\alpha([0, T], \xR^n) \subset \C([0, T], \xR^n)$ and $\Lip([0, T], \xR) \subset \C([0, T], \xR)$, the set $XZ^\alpha(t_0, w_0(\cdot), z_0)$ can be considered as a subset of $\C([0, T], \xR^n \times \xR)$.

        By a solution of problem \eqref{Differential_Inclusion}, \eqref{Differential_Inclusion_Initial_Condition}, we mean a pair of functions $(x(\cdot), z(\cdot)) \in XZ^\alpha(t_0, w_0(\cdot), z_0)$ such that differential inclusion \eqref{Differential_Inclusion} is fulfilled for a.e. $t \in [t_0, T]$.
        Let $XZ^\alpha_0(t_0, w_0(\cdot), z_0)$ denote the set of such solutions.

        \begin{prpstn} \label{proposition_DI_1}
            For any $(t_0, w_0(\cdot)) \in G_n^\alpha$ and $z_0 \in \xR$, the set $XZ^\alpha_0(t_0, w_0(\cdot), z_0)$ is nonempty and compact in $\C([0, T], \xR^n \times \xR)$.
        \end{prpstn}
        \begin{prpstn} \label{proposition_DI_2}
            Let $(t_k, w_k(\cdot)) \in G_n^\alpha$, $z_k \in \xR$, and $(x^{(k)}(\cdot), z^{(k)}(\cdot)) \in XZ^\alpha_0(t_k, w_k(\cdot), z_k)$ for every $k \in \xN$, and let $(t_k, w_k(\cdot)) \to (t_0, w_0(\cdot)) \in G_n^\alpha$ and $z_k \to z_0 \in \xR$ as $k \to \infty$.
            Then, the sequence $\{(x^{(k)}(\cdot), z^{(k)}(\cdot))\}_{k \in \xN}$ contains a subsequence that converges to a solution $(x^{(0)}(\cdot), z^{(0)}(\cdot)) \in XZ^\alpha_0(t_0, w_0(\cdot), z_0)$.
        \end{prpstn}
        \begin{prpstn} \label{proposition_semigroup_property}
            Let $(t_0, w_0(\cdot)) \in G_n^\alpha$, $z_0 \in \xR$, and $(x(\cdot), z(\cdot)) \in XZ^\alpha_0(t_0, w_0(\cdot), z_0)$.
            Then, for every $t^\prime \in [t_0, T]$ and $(x^\prime(\cdot), z^\prime(\cdot)) \in XZ^\alpha_0(t^\prime, x_{t^\prime}(\cdot), z(t^\prime))$, the inclusion $(x^\prime(\cdot), z^{\prime \prime}(\cdot)) \in XZ^\alpha_0(t_0, w_0(\cdot), z_0)$ holds, where the function $z^{\prime \prime}(\cdot)$ is defined by $z^{\prime \prime}(t) = z(t)$ for $t \in [0, t^\prime]$ and $z^{\prime \prime}(t) = z^\prime(t)$ for $t \in (t^\prime, T]$.
        \end{prpstn}

        In the case when there is no additional variable $z(t)$, similar statements are proved in \cite{Gomoyunov_2020_DE} by suitably adapting the proofs of the corresponding results for ordinary and functional differential inclusions with first-order derivatives (see, e.g., \cite{Filippov_1988} and also \cite{Kurzhanskii_1970_Eng,Lukoyanov_2011_Eng}).
        The proofs of Propositions~\ref{proposition_DI_1},~\ref{proposition_DI_2}, and~\ref{proposition_semigroup_property} can be carried out by the same scheme with only minor technical changes, and, therefore, they are omitted.

        Finally, note that, since the right-hand side of differential inclusion \eqref{Differential_Inclusion} does not depend on $z(t)$, then, for any $(t_0, w_0(\cdot)) \in G_n^\alpha$, $z_0 \in \xR$, $(x(\cdot), z(\cdot)) \in XZ^\alpha_0(t_0, w_0(\cdot), z_0)$, and $z^\prime_0 \in \xR$, we have $(x(\cdot), z^\prime(\cdot)) \in XZ^\alpha_0(t_0, w_0(\cdot), z_0^\prime)$ for the function $z^\prime(t) = z^\prime_0 + z(t) - z_0$, $t \in [0, T]$.

    \subsection{Functional Differential Inclusions}

        Let us also give an analogue of Proposition~\ref{proposition_DI_1} for the case when the right-hand side of differential inclusion \eqref{Differential_Inclusion} depends not only on a single value $x(t)$ of an unknown solution, but also on all values $x(\tau)$, $\tau \in [0, t]$, or, in other words, on the function $x_t(\cdot)$ given by \eqref{x_t}.

        Let a set-valued functional $G_n^\alpha \ni (t, w(\cdot)) \mapsto \mathcal{F}(t, w(\cdot)) \subset \xR^n \times \xR$ be such that:
        \begin{itemize}
            \item[$(\mathcal{F}.1)$]
                For every $(t, w(\cdot)) \in G_n^\alpha$, the set $\mathcal{F}(t, w(\cdot))$ is nonempty, convex, and compact in $\xR^n \times \xR$.

            \item[$(\mathcal{F}.2)$]
                The set-valued functional $\mathcal{F}$ is upper semicontinuous.
                Namely, for every $(t, w(\cdot)) \in G_n^\alpha$ and $\varepsilon > 0$, there exists $\delta > 0$ such that, for any $(t^\prime, w^\prime(\cdot)) \in G_n^\alpha$, the inequality $\dist((t, w(\cdot)), (t^\prime, w^\prime(\cdot))) \leq \delta$ implies the inclusion $\mathcal{F}(t^\prime, w^\prime(\cdot)) \subset [\mathcal{F}(t, w(\cdot))]^\varepsilon$.

            \item[$(\mathcal{F}.3)$]
                There exists $c_{\mathcal{F}} > 0$ such that, for any $(t, w(\cdot)) \in G_n^\alpha$,
                \begin{equation*}
                    \sup \big\{\|f\|:
                    \, (f, h) \in \mathcal{F}(t, w(\cdot))\big\}
                    \leq c_F (1 + \|w(\cdot)\|_{[0, t]}).
                \end{equation*}
        \end{itemize}

        Given $(t_0, w_0(\cdot)) \in G_n^\alpha$ and $z_0 \in \xR$, consider a Cauchy problem for the functional differential inclusion
        \begin{equation} \label{Functional_Differential_Inclusion}
            \big( (^C D^\alpha x)(t), \dot{z}(t) \big)
            \in \mathcal{F}(t, x_t(\cdot)),
            \quad (x(t), z(t)) \in \xR^n \times \xR, \quad t \in [t_0, T],
        \end{equation}
        and initial condition \eqref{Differential_Inclusion_Initial_Condition}.
        A pair of functions $(x(\cdot), z(\cdot)) \in XZ^\alpha(t_0, w_0(\cdot), z_0)$ is a solution of this problem if functional differential inclusion \eqref{Functional_Differential_Inclusion} holds for a.e. $t \in [t_0, T]$.
        Let $\mathcal{XZ}^\alpha_0(t_0, w_0(\cdot), z_0)$ be the set of such solutions.

        \begin{prpstn} \label{proposition_FDI}
            For any $(t_0, w_0(\cdot)) \in G_n^\alpha$ and $z_0 \in \xR$, the set $\mathcal{XZ}^\alpha_0(t_0, w_0(\cdot), z_0)$ is nonempty and compact in $\C([0, T], \xR^n \times \xR)$.
        \end{prpstn}

        This proposition can be proved by the scheme from \cite[Theorem~1]{Gomoyunov_2020_DE}.

\section{Hamilton--Jacobi Equation with Fractional Coinvariant Derivatives}
\label{section_HJ}

    In this section, we consider a Cauchy problem for a Hamilton--Jacobi equation with fractional coinvariant ($ci$-) derivatives of an order $\alpha \in (0, 1)$ and propose a definition of a minimax solution of this problem.

    \subsection{Fractional Coinvariant Derivatives}
    \label{subsection_ci-derivatives}

        Let us recall the notion of $ci$-differentiability of the order $\alpha$ of a functional $\varphi: G^\alpha_n \to \xR$, which was introduced in \cite{Gomoyunov_2020_SIAM} as a suitable extension of the notion of $ci$-differentiability (of the first order) developed in, e.g., \cite{Kim_1999,Lukoyanov_2011_Eng}.

        For a given $(t_0, w_0(\cdot)) \in G_n^{\alpha \circ}$, consider the set of admissible extensions $x(\cdot)$ of $w_0(\cdot)$ defined as follows:
        \begin{equation} \label{X^alpha}
            X^\alpha(t_0, w_0(\cdot))
            = \big\{ x(\cdot) \in \AC^\alpha([0, T], \xR^n):
            \, x(t) = w_0(t), \, t \in [0, t_0] \big\}.
        \end{equation}
        The functional $\varphi$ is called $ci$-differentiable of the order $\alpha$ at $(t_0, w_0(\cdot))$ if there exist $\partial_t^\alpha \varphi(t_0, w_0(\cdot)) \in \xR$ and $\nabla^\alpha \varphi(t_0, w_0(\cdot)) \in \xR^n$ such that, for every extension $x(\cdot) \in X^\alpha(t_0, w_0(\cdot))$, the relation
        \begin{equation} \label{ci_differentiability}
            \varphi(t, x_t(\cdot)) - \varphi(t_0, w_0(\cdot))
            = \partial_t^\alpha \varphi(t_0, w_0(\cdot)) (t - t_0)
            + \int_{t_0}^{t} \langle \nabla^\alpha \varphi(t_0, w_0(\cdot)), (^C D^\alpha x)(\tau) \rangle \xdif \tau
            + o(t - t_0)
        \end{equation}
        holds for any $t \in (t_0, T)$.
        Here, $x_t(\cdot)$ is determined by $x(\cdot)$ and $t$ according to \eqref{x_t}, the function $o$ may depend on $t$ and $x(\cdot)$, and $o(\delta) / \delta \to 0$ as $\delta \to + 0$.
        In this case, the quantities $\partial_t^\alpha \varphi(t_0, w_0(\cdot))$ and $\nabla^\alpha \varphi(t_0, w_0(\cdot))$ are called the $ci$-derivatives of the order $\alpha$ of $\varphi$ at $(t_0, w_0(\cdot))$.

        The functional $\varphi$ is said to be $ci$-smooth of the order $\alpha$ if it is continuous, $ci$-differentiable of the order $\alpha$ at every $(t, w(\cdot)) \in G_n^{\alpha \circ}$, and the functionals $\partial_t^\alpha \varphi: G_n^{\alpha \circ} \to \xR$ and $\nabla^\alpha \varphi: G_n^{\alpha \circ} \to \xR^n$ are continuous.
        Recall that the set $G_n^{\alpha \circ} \subset G_n$ is endowed with the metric $\dist$ from \eqref{dist}.

    \subsection{Hamilton--Jacobi Equation}
    \label{subsection_HJE}

        Consider the Cauchy problem for the Hamilton--Jacobi equation with $ci$-derivatives of the order $\alpha$
        \begin{equation} \label{HJ}
            \partial^\alpha_t \varphi(t, w(\cdot))
            + H\big( t, w(t), \nabla^\alpha \varphi(t, w(\cdot)) \big)
            = 0,
            \quad (t, w(\cdot)) \in G_n^{\alpha \circ},
        \end{equation}
        and the boundary condition
        \begin{equation}\label{HJ_boundary_condition}
            \varphi(T, w(\cdot))
            = \sigma(w(\cdot)),
            \quad w(\cdot) \in \AC^\alpha([0, T], \xR^n).
        \end{equation}
        In this problem, $\varphi: G_n^\alpha \to \xR$ is an unknown functional, and the given mappings $H: [0, T] \times \xR^n \times \xR^n \to \xR$ and $\sigma: \AC^\alpha([0, T], \xR^n) \to \xR$ are assumed to satisfy the following conditions:
        \begin{itemize}
            \item[$(H.1)$]
                The function $H$ is continuous.

            \item[$(H.2)$]
                There exists $c_H > 0$ such that, for any $t \in [0, T]$ and $x$, $s$, $s^\prime \in \xR^n$,
                \begin{equation*}
                    |H(t, x, s) - H(t, x, s^\prime)|
                    \leq c_H (1 + \|x\|) \|s - s^\prime\|.
                \end{equation*}

            \item[$(H.3)$]
                For every $R \geq 0$, there exists $\lambda_H > 0$ such that, for any $t \in [0, T]$ and $x$, $x^\prime$, $s \in \xR^n$, if $\|x\| \leq R$ and $\|x^\prime\| \leq R$, then
                \begin{equation*}
                    |H(t, x, s) - H(t, x^\prime, s)|
                    \leq \lambda_H (1 + \|s\|) \|x - x^\prime\|.
                \end{equation*}

            \item[$(\sigma)$]
                 The functional $\sigma$ is continuous.
        \end{itemize}

        Cauchy problem \eqref{HJ}, \eqref{HJ_boundary_condition} arises \cite{Gomoyunov_2020_SIAM} when studying infinitesimal properties of the value functional in Bolza-type optimal control problems for dynamical systems described by fractional differential equations with the Caputo derivatives of the order $\alpha$.
        In this connection, assumptions $(H.1)$--$(H.3)$ and $(\sigma)$ seem quite natural since they are fulfilled in a sufficiently wide range of such problems.
        If the value functional is $ci$-smooth of the order $\alpha$, then, according to \cite[Theorem~10.1]{Gomoyunov_2020_SIAM}, it satisfies Hamilton--Jacobi equation \eqref{HJ} and boundary condition \eqref{HJ_boundary_condition}, and, therefore, it can be considered as a solution of problem \eqref{HJ}, \eqref{HJ_boundary_condition} in the classical sense.
        In particular, this allows us to efficiently construct optimal control strategies \cite[Corollary~11.4]{Gomoyunov_2020_SIAM}.
        However, the value functional usually does not possess such smoothness properties, which leads to the need to introduce and study generalized solutions of problem \eqref{HJ}, \eqref{HJ_boundary_condition}.

        Below, we give a definition of a minimax solution of problem \eqref{HJ}, \eqref{HJ_boundary_condition}, which is a suitable modification of the corresponding definitions is the case of Hamilton--Jacobi equations with partial derivatives (see, e.g., \cite[Sect.~6.2]{Subbotin_1995}) and with first-order $ci$-derivatives (see, e.g., \cite[Sect.~6]{Lukoyanov_2011_Eng}, \cite{Lukoyanov_2003_1}, and also \cite{Bayraktar_Keller_2018}).
        We prove that the minimax solution exists, is unique, and is consistent with a classical solution of problem \eqref{HJ}, \eqref{HJ_boundary_condition}.
        In particular, we establish a comparison principle.

    \subsection{Minimax Solution}
    \label{subsection_minimax_solution}

        Consider the set-valued function $[0, T] \times \xR^n \times \xR^n \ni (t, x, s) \mapsto E(t, x, s) \subset \xR^n \times \xR$, where
        \begin{equation} \label{E}
            E(t, x, s)
            = \big\{ (f, h) \in \xR^n \times \xR:
            \, \|f\| \leq c_H(1 + \|x\|),
            \, h = \langle s, f \rangle - H(t, x, s) \big\},
            \quad t \in [0, T], \quad x, s \in \xR^n.
        \end{equation}
        Note that (see, e.g., \cite[Sect.~6.2]{Subbotin_1995}) the set $E(t, x, s)$ is nonempty, convex, and compact in $\xR^n \times \xR$ for every $t \in [0, T]$ and $x$, $s \in \xR^n$, the set-valued function $E$ is continuous (in the Hausdorff sense) due to assumption $(H.1)$, and the inequality below holds:
        \begin{equation*}
            \sup \big\{\|f\|:
            \, (f, h) \in E(t, x, s)\big\}
            \leq c_H (1 + \|x\|),
            \quad t \in [0, T], \quad x, s \in \xR^n.
        \end{equation*}
        In addition, it follows from $(H.2)$ that $E(t, x, s) \cap E(t, x, s^\prime) \neq \emptyset$ for any $t \in [0, T]$ and $x$, $s$, $s^\prime \in \xR^n$.

        Given $(t_0, w_0(\cdot)) \in G_n^\alpha$, $z_0 \in \xR$, and $s \in \xR^n$, consider the Cauchy problem for the differential inclusion
        \begin{equation} \label{CH_DI}
            \big( (^C D^\alpha x)(t), \dot{z}(t) \big)
            \in E(t, x(t), s),
            \quad (x(t), z(t)) \in \xR^n \times \xR, \quad t \in [t_0, T],
        \end{equation}
        and the initial condition
        \begin{equation} \label{CH_DI_initial_condition}
            x(t)
            = w_0(t),
            \quad z(t)
            = z_0,
            \quad t \in [0, t_0].
        \end{equation}
        Here, $s$ is treated as a constant parameter.
        Let $CH(t_0, w_0(\cdot), z_0, s)$ be the set of solutions $(x(\cdot), z(\cdot))$ of problem \eqref{CH_DI}, \eqref{CH_DI_initial_condition}.
        According to Proposition~\ref{proposition_DI_1} and the described above properties of the function $E$, the set $CH(t_0, w_0(\cdot), z_0, s)$ is nonempty and compact in $\C([0, T], \xR^n \times \xR)$.
        Following the conventional terminology, differential inclusion \eqref{CH_DI} is called a characteristic differential inclusion, and any element of the set $CH(t_0, w_0(\cdot), z_0, s)$ is called a (generalized) characteristic of equation \eqref{HJ}.

        We say that a functional $\varphi: G_n^\alpha \to \xR$ is an upper solution of problem \eqref{HJ}, \eqref{HJ_boundary_condition} if it is lower semicontinuous, satisfies the boundary condition
        \begin{equation} \label{upper_solution_boundary}
            \varphi(T, w(\cdot))
            \geq \sigma(w(\cdot)),
            \quad w(\cdot) \in \AC^\alpha([0, T], \xR^n),
        \end{equation}
        and possesses the following property:
        \begin{itemize}
            \item[$(\varphi_+)$]
                For every $(t_0, w_0(\cdot)) \in G_n^{\alpha \circ}$, $t \in (t_0, T]$, $s \in \xR^n$, and $\varepsilon > 0$, there exists $(x(\cdot), z(\cdot)) \in CH(t_0, w_0(\cdot), 0, s)$ such that $\varphi(t, x_t(\cdot)) - z(t) \leq \varphi(t_0, w_0(\cdot)) + \varepsilon$.
        \end{itemize}
        Respectively, a lower solution of this problem is an upper semicontinuous functional $\varphi: G_n^\alpha \to \xR$ such that
        \begin{equation} \label{lower_solution_boundary}
            \varphi(T, w(\cdot))
            \leq \sigma(w(\cdot)),
            \quad w(\cdot) \in \AC^\alpha([0, T], \xR^n),
        \end{equation}
        and the statement below holds:
        \begin{itemize}
            \item[$(\varphi_-)$]
                For every $(t_0, w_0(\cdot)) \in G_n^{\alpha \circ}$, $t \in (t_0, T]$, $s \in \xR^n$, and $\varepsilon > 0$, there exists $(x(\cdot), z(\cdot)) \in CH(t_0, w_0(\cdot), 0, s)$ such that $\varphi(t, x_t(\cdot)) - z(t) \geq \varphi(t_0, w_0(\cdot)) - \varepsilon$.
        \end{itemize}
        In $(\varphi_+)$ and $(\varphi_-)$, as usual, the function $x_t(\cdot)$ is the restriction of the function $x(\cdot)$ to the interval $[0, t]$ (see \eqref{x_t}).

        A functional $\varphi: G_n^\alpha \to \xR$ is called a minimax solution of problem \eqref{HJ}, \eqref{HJ_boundary_condition} if it is an upper solution as well as a lower solution of this problem.

        \begin{rmrk}
            Conditions $(\varphi_+)$ and $(\varphi_-)$ can be reformulated in terms of weak invariance of respectively the epigraph and hypograph of the functional $\varphi$ with respect to characteristic differential inclusion \eqref{CH_DI} for every $s \in \xR^n$ (see, e.g., definitions $(U2)$ and $(L2)$ in \cite[Sect~6.3]{Subbotin_1995}).
            Note also that, in the terminology of positional differential games theory, statements $(\varphi_+)$ and $(\varphi_-)$ express so-called $u$-stability and $v$-stability properties of the value function (see, e.g., \cite[Sect.~4.2]{Krasovskii_Subbotin_1988} and \cite[Sect.~8]{Krasovskii_Krasovskii_1995}).
        \end{rmrk}

\section{Consistency}
    \label{section_Consistency}

    This section deals with issues of consistency of a minimax solution of problem \eqref{HJ}, \eqref{HJ_boundary_condition} with a solution of this problem in the classical sense.

    By a classical solution of problem \eqref{HJ}, \eqref{HJ_boundary_condition}, we mean a $ci$-smooth of the order $\alpha$ functional $\varphi: G_n^\alpha \to \xR$ that satisfies Hamilton--Jacobi equation \eqref{HJ} and boundary condition \eqref{HJ_boundary_condition}.

    The schemes of the proof of the statements below go back to the proofs of the corresponding results for Hamilton--Jacobi equations with partial derivatives (see, e.g., \cite[Sect.~2.4]{Subbotin_1995}) and with first-order $ci$-derivatives (see, e.g., \cite[Sect.~4 and~5]{Lukoyanov_2011_Eng}, \cite[Proposition~5.1]{Lukoyanov_2003_1}, and also \cite[Sect.~B.1]{Bayraktar_Keller_2018}).

    \begin{thrm} \label{theorem_classical_is_minimax}
        A classical solution of problem \eqref{HJ}, \eqref{HJ_boundary_condition} is a minimax solution of this problem.
    \end{thrm}
    \begin{proof}
        Since a classical solution $\varphi: G_n^\alpha \to \xR$ is continuous and satisfies \eqref{HJ_boundary_condition}, in order to prove that $\varphi$ is a minimax solution, it suffices to verify that $\varphi$ possesses properties $(\varphi_+)$ and $(\varphi_-)$.
        To this end, let us show that, for a given $(t_0, w_0(\cdot)) \in G_n^{\alpha \circ}$ and $s \in \xR^n$, there exists $(x^\ast(\cdot), z^\ast(\cdot)) \in CH(t_0, w_0(\cdot), 0, s)$ such that
        \begin{equation} \label{theorem_classical_is_minimax_main}
            \varphi(t, x^\ast_t(\cdot)) - z^\ast(t)
            = \varphi(t_0, w_0(\cdot)),
            \quad t \in [t_0, T].
        \end{equation}
        Consider the set-valued functional $G_n^\alpha \ni (t, w(\cdot)) \mapsto \mathcal{E}^\ast(t, w(\cdot)) \subset \xR^n \times \xR$, where, for every $(t, w(\cdot)) \in G_n^\alpha$,
        \begin{equation*}
            \mathcal{E}^\ast(t, w(\cdot))
            = \begin{cases}
                E(t, w(t), s) \cap E\big(t, w(t), \nabla^\alpha \varphi(t, w(\cdot))\big), & \mbox{if } t < T, \\
                E(t, w(t), s), & \mbox{if } t = T.
              \end{cases}
        \end{equation*}
        Due to the given in Sect.~\ref{subsection_minimax_solution} properties of the set-valued function $E$ and continuity of the functionals $\nabla^\alpha \varphi$ and $G_n \ni (t, w(\cdot)) \mapsto w(t) \in \xR^n$, the functional $\mathcal{E}^\ast$ satisfies conditions $(\mathcal{F}.1)$--$(\mathcal{F}.3)$.
        Then, owing to Proposition~\ref{proposition_FDI}, the Cauchy problem for the functional differential inclusion
        \begin{equation*}
            \big( (^C D^\alpha x)(t), \dot{z}(t) \big)
            \in \mathcal{E}^\ast(t, x_t(\cdot)),
            \quad (x(t), z(t)) \in \xR^n \times \xR, \quad t \in [t_0, T],
        \end{equation*}
        and the initial condition $x(t) = w_0(t)$, $z(t) = 0$, $t \in [0, t_0]$, admits a solution $(x^\ast(\cdot), z^\ast(\cdot))$.
        By construction, we have $(x^\ast(\cdot), z^\ast(\cdot)) \in CH(t_0, w_0(\cdot), 0, s)$ and, in accordance with \eqref{E},
        \begin{equation*}
            \dot{z}^\ast(t)
            = \langle \nabla^\alpha \varphi(t, x^\ast_t(\cdot)), (^C D^\alpha x^\ast)(t) \rangle
            - H\big(t, x^\ast(t), \nabla^\alpha \varphi(t, x^\ast_t(\cdot))\big)
            \text{ for a.e. } t \in [t_0, T].
        \end{equation*}
        Hence, taking into account that $\varphi$ satisfies equation \eqref{HJ}, we get
        \begin{equation*}
            \dot{z}^\ast(t)
            = \partial^\alpha_t \varphi(t, x^\ast_t(\cdot)) + \langle \nabla^\alpha \varphi(t, x^\ast_t(\cdot)), (^C D^\alpha x^\ast)(t) \rangle
            \text{ for a.e. } t \in [t_0, T].
        \end{equation*}
        On the other hand, since $\varphi$ is $ci$-smooth of the order $\alpha$, by \cite[Lemma~9.2]{Gomoyunov_2020_SIAM}, for the function $\omega(t) = \varphi(t, x^\ast_t(\cdot))$, $t \in [t_0, T]$, and a fixed $\vartheta \in [t_0, T)$, we obtain
        \begin{equation*}
            \omega(\vartheta)
            = \omega(t_0) + \int_{t_0}^{\vartheta} \big( \partial^\alpha_t \varphi(t, x^\ast_t(\cdot))
            + \langle \nabla^\alpha \varphi(t, x^\ast_t(\cdot)), (^C D^\alpha x^\ast)(t) \rangle \big) \xdif t.
        \end{equation*}
        Thus, recalling that $x^\ast_{t_0}(\cdot) = w_0(\cdot)$, $z^\ast(t_0) = 0$, and $z^\ast(\cdot) \in \Lip([0, T], \xR)$, we derive
        \begin{equation*}
            \varphi(\vartheta, x^\ast_\vartheta(\cdot))
            = \varphi(t_0, w_0(\cdot)) + \int_{t_0}^{\vartheta} \dot{z}^\ast(t) \xdif t
            = \varphi(t_0, w_0(\cdot)) + z^\ast(\vartheta).
        \end{equation*}
        This equality is valid for every $\vartheta \in [t_0, T)$, and, therefore, in view of continuity of $\varphi$, we conclude \eqref{theorem_classical_is_minimax_main}, which completes the proof of the theorem.
    \end{proof}

    We also establish the following result.
    \begin{thrm} \label{theorem_minimax_is_differentiable}
        If a minimax solution of problem \eqref{HJ}, \eqref{HJ_boundary_condition} is $ci$-differentiable of the order $\alpha$ at some point $(t_0, w_0(\cdot)) \in G_n^{\alpha \circ}$, then it satisfies equation \eqref{HJ} at this point.
    \end{thrm}

    Before proving the theorem, we present an auxiliary proposition.
    \begin{prpstn} \label{proposition_strong_stability}
        If a functional $\varphi: G_n^\alpha \to \xR$ is lower semicontinuous, then $\varphi$ satisfies condition $(\varphi_+)$ if and only if the following statement holds:
        \begin{itemize}
            \item[$(\varphi_+^\ast)$]
                For any $(t_0, w_0(\cdot)) \in G_n^{\alpha \circ}$ and $s \in \xR^n$, there is a characteristic $(x(\cdot), z(\cdot)) \in CH(t_0, w_0(\cdot), 0, s)$ such that $\varphi(t, x_t(\cdot)) - z(t) \leq \varphi(t_0, w_0(\cdot))$ for every $t \in [t_0, T]$.
        \end{itemize}
        Respectively, for an upper semicontinuous functional $\varphi: G_n^\alpha \to \xR$, condition $(\varphi_-)$ is equivalent to the following:
        \begin{itemize}
            \item[$(\varphi_-^\ast)$]
                For any $(t_0, w_0(\cdot)) \in G_n^{\alpha \circ}$ and $s \in \xR^n$, there is a characteristic $(x(\cdot), z(\cdot)) \in CH(t_0, w_0(\cdot), 0, s)$ such that $\varphi(t, x_t(\cdot)) - z(t) \geq \varphi(t_0, w_0(\cdot))$ for every $t \in [t_0, T]$.
        \end{itemize}
    \end{prpstn}
    \begin{proof}
        We prove only the first part of the proposition since the proof of the second one is essentially the same.
        It is clear that $(\varphi_+)$ follows from $(\varphi_+^\ast)$, so it remains to verify the reverse implication.
        Let $\varphi: G_n^\alpha \to \xR$ be a lower semicontinuous functional satisfying $(\varphi_+)$, and let $(t_0, w_0(\cdot)) \in G_n^{\alpha \circ}$ and $s \in \xR^n$.

        Fix $k \in \xN$.
        Denote $t_{k, i} = t_0 + (T - t_0) i / k$, $i \in \overline{0, k}$.
        Take arbitrarily $(x^{(k, 0)}(\cdot), z^{(k, 0)}(\cdot)) \in CH(t_0, w_0(\cdot), 0, s)$ and, applying $(\varphi_+)$, choose functions $x^{(k, i)}: [0, T] \to \xR^n$, $z^{(k, i)}: [0, T] \to \xR$, $i \in \overline{1, k}$, such that the following relations hold for every $i \in \overline{1, k}$:
        \begin{equation*}
            (x^{(k, i)}(\cdot), z^{(k, i)}(\cdot))
            \in CH\big(t_{k, i - 1}, x^{(k, i - 1)}_{t_{k, i - 1}}(\cdot), 0, s\big),
            \quad \varphi\big(t_{k, i}, x^{(k, i)}_{t_{k, i}}(\cdot)\big) - z^{(k, i)}(t_{k, i})
            \leq \varphi\big(t_{k, i - 1}, x^{(k, i - 1)}_{t_{k, i - 1}}(\cdot)\big) + 1 / k^2.
        \end{equation*}
        Further, consider functions $\bar{z}^{(k, i)}: [0, T] \to \xR$, $i \in \overline{0, k}$, such that $\bar{z}^{(k, 0)}(\cdot) = z^{(k, 0)}(\cdot)$ and, for any $i \in \overline{1, k}$,
        \begin{equation*}
            \bar{z}^{(k, i)}(t)
            = \begin{cases}
                \bar{z}^{(k, i - 1)}(t), & \mbox{if } t \in [0, t_{k, i - 1}], \\
                z^{(k, i)}(t) + \bar{z}^{(k, i - 1)}(t_{k, i - 1}), & \mbox{if } t \in (t_{k, i - 1}, T].
              \end{cases}
        \end{equation*}
        Then, by induction, based on Proposition~\ref{proposition_semigroup_property} (see also the remark in the end of Sect.~\ref{subsection_ODI}), we can prove that, for every $i \in \overline{0, k}$,
        \begin{equation*}
            (x^{(k, i)}(\cdot), \bar{z}^{(k, i)}(\cdot))
            \in CH(t_0, w_0(\cdot), 0, s),
            \quad \varphi\big(t_{k, j}, x^{(k, i)}_{t_{k, j}}(\cdot)\big) - \bar{z}^{(k, i)}(t_{k, j})
            \leq \varphi(t_0, w_0(\cdot)) + j / k^2,
            \quad j \in \overline{0, i}.
        \end{equation*}
        Thus, for the functions $x^{[k]}(\cdot) = x^{(k, k)}(\cdot)$ and $z^{[k]}(\cdot) = \bar{z}^{(k, k)}(\cdot)$, we obtain
        \begin{equation*}
            (x^{[k]}(\cdot), z^{[k]}(\cdot))
            \in CH(t_0, w_0(\cdot), 0, s),
            \quad \varphi\big(t_{k, j}, x^{[k]}_{t_{k, j}}(\cdot)\big) - z^{[k]}(t_{k, j})
            \leq \varphi(t_0, w_0(\cdot)) + 1 / k,
            \quad j \in \overline{0, k}.
        \end{equation*}

        Due to compactness of $CH(t_0, w_0(\cdot), 0, s)$, we can assume that the sequence $\{(x^{[k]}(\cdot), z^{[k]}(\cdot))\}_{k \in \xN}$ converges to a characteristic $(x^{[0]}(\cdot), z^{[0]}(\cdot)) \in CH(t_0, w_0(\cdot), 0, s)$.
        Now, let $t \in [t_0, T]$ be fixed.
        For every $k \in \xN$, denoting $t^{[k]} = \max \{ t_{k, i}: \, t_{k, i} \leq t, \, i \in \overline{0, k}\}$, we get
        \begin{equation} \label{proposition_strong_stability_proof_main}
            \varphi\big(t^{[k]}, x^{[k]}_{t^{[k]}}(\cdot)\big) - z^{[k]}(t^{[k]})
            \leq \varphi(t_0, w_0(\cdot)) + 1 / k.
        \end{equation}
        As $k \to \infty$, we have $t^{[k]} \to t$, $(t^{[k]}, x^{[k]}_{t^{[k]}}(\cdot)) \to (t, x^{[0]}_t(\cdot))$, and $z^{[k]}(t^{[k]}) \to z^{[0]}(t)$.
        Hence, passing to the limit as $k \to \infty$ in inequality \eqref{proposition_strong_stability_proof_main}, in view of lower semicontinuity of $\varphi$, we derive
        \begin{equation*}
            \varphi(t, x^{[0]}_t(\cdot)) - z^{[0]}(t)
            \leq \liminf_{k \to \infty} \big( \varphi\big(t^{[k]}, x^{[k]}_{t^{[k]}}(\cdot)\big) - z^{[k]}(t^{[k]}) \big)
            \leq \varphi(t_0, w_0(\cdot)).
        \end{equation*}
        So, the functional $\varphi$ possesses property $(\varphi_+^\ast)$, and the proposition is proved.
    \end{proof}

    \begin{proof}[Proof of Theorem~\ref{theorem_minimax_is_differentiable}.]
        Assume that $(t_0, w_0(\cdot)) \in G_n^{\alpha \circ}$ and a minimax solution $\varphi: G_n^\alpha \to \xR$ of problem \eqref{HJ}, \eqref{HJ_boundary_condition} is $ci$-differentiable of the order $\alpha$ at $(t_0, w_0(\cdot))$.
        Denote $s_0 = \nabla^\alpha \varphi(t_0, w_0(\cdot))$.
        Since $\varphi$ is an upper solution of problem \eqref{HJ}, \eqref{HJ_boundary_condition}, due to $(\varphi_+^\ast)$, there exists a characteristic $(x(\cdot), z(\cdot)) \in CH(t_0, w_0(\cdot), 0, s_0)$ such that $\varphi (t, x_t(\cdot)) - z(t) \leq \varphi (t_0, w_0(\cdot))$, $t \in [t_0, T]$.
        In particular, according to \eqref{E}, we have
        \begin{equation*}
            z(t)
            = \int_{t_0}^{t} \big( \langle s_0, (^C D^\alpha x)(\tau) \rangle - H(\tau, x(\tau), s_0) \big) \xdif \tau,
            \quad t \in [t_0, T].
        \end{equation*}
        Hence, taking into account that $x(\cdot) \in X^\alpha(t_0, w_0(\cdot))$, in view of \eqref{ci_differentiability}, we derive
        \begin{equation} \label{theorem_minimax_is_differentiable_proof_main}
            0
            \geq \varphi (t, x_t(\cdot)) - z(t) - \varphi (t_0, w_0(\cdot))
            = \partial^\alpha_t \varphi(t_0, w_0(\cdot)) (t - t_0) + \int_{t_0}^{t} H(\tau, x(\tau), s_0) \xdif \tau
            + o(t - t_0),
            \quad t \in (t_0, T),
        \end{equation}
        Note that $H(t, x(t), s_0) \to H(t_0, w_0(t_0), s_0)$ as $t \to t_0 + 0$ by virtue of assumption $(H.1)$.
        Therefore, dividing \eqref{theorem_minimax_is_differentiable_proof_main} by $t - t_0$ and, after that, passing to the limit as $t \to t_0 + 0$, we get
        \begin{equation} \label{theorem_minimax_is_differentiable_proof_first}
            0
            \geq \partial^\alpha_t \varphi(t_0, w_0(\cdot)) + H(t_0, w_0(t_0), s_0).
        \end{equation}
        On the other hand, based on the fact that $\varphi$ is a lower solution of problem \eqref{HJ}, \eqref{HJ_boundary_condition}, and, consequently, it possesses property $(\varphi_-^\ast)$, we can similarly obtain the inequality
        \begin{equation} \label{theorem_minimax_is_differentiable_proof_second}
            0
            \leq \partial^\alpha_t \varphi(t_0, w_0(\cdot)) + H(t_0, w_0(t_0), s_0).
        \end{equation}
        It follows from \eqref{theorem_minimax_is_differentiable_proof_first} and \eqref{theorem_minimax_is_differentiable_proof_second} that $\varphi$ satisfies equation \eqref{HJ} at $(t_0, w_0(\cdot))$.
        The theorem is proved.
    \end{proof}

    In particular, from Theorem~\ref{theorem_minimax_is_differentiable}, we derive
    \begin{crllr} \label{corollary_minimax_is_classical}
        If a minimax solution of problem \eqref{HJ}, \eqref{HJ_boundary_condition} is $ci$-smooth of the order $\alpha$, then it is a classical solution of this problem.
    \end{crllr}

    Theorems~\ref{theorem_classical_is_minimax} and~\ref{theorem_minimax_is_differentiable} and Corollary~\ref{corollary_minimax_is_classical} allow us to conclude that the introduced notion of a minimax solution of problem \eqref{HJ}, \eqref{HJ_boundary_condition} is consistent with the notion of a solution of this problem in the classical sense.

\section{Comparison Principle}
\label{section_Comparison}

    The goal of this section is to prove the result below, which is often called a comparison principle.
    In the next section, it is used in the proof of existence and uniqueness of a minimax solution of problem \eqref{HJ}, \eqref{HJ_boundary_condition}.
    \begin{thrm} \label{theorem_comparison_priciple}
        Let $\varphi_+$ and $\varphi_-$ be respectively an upper and a lower solutions of problem \eqref{HJ}, \eqref{HJ_boundary_condition}.
        Then, the inequality below holds:
        \begin{equation} \label{theorem_comparison_priciple_main}
            \varphi_-(t, w(\cdot))
            \leq \varphi_+(t, w(\cdot)),
            \quad (t, w(\cdot)) \in G_n^\alpha.
        \end{equation}
    \end{thrm}

    In general, this theorem is proved by the same scheme as the corresponding statements for Hamilton--Jacobi equations with partial derivatives (see, e.g., \cite[Theorem~7.3]{Subbotin_1995}) and with first-order $ci$-derivatives (see, e.g., \cite[Lemma~7.7]{Lukoyanov_2003_1}).
    However, the key point of the proof, which concerns construction of a Lyapunov--Krasovskii functional with a number of prescribed properties (in this connection, see, e.g., \cite[Sect.~5]{Lukoyanov_2004_PMM_Eng}), substantially differs from the previous studies owing to features of fractional-order integrals and derivatives.

    \subsection{Lyapunov--Krasovskii Functionals}

        The construction of the required Lyapunov--Krasovskii functional is carried out in four steps.

        \subsubsection{Functional $V_{\gamma, \mu}$}

        Given $\gamma \in (0, 1)$ and $\mu > 0$, consider the functional
        \begin{equation} \label{V^1_gamma_mu}
            G_1 \ni (t, r(\cdot))
            \mapsto V_{\gamma, \mu}(t, r(\cdot))
            = \frac{1}{\Gamma(1 - \gamma)} \int_{0}^{t} \frac{e^{- \mu (t - \tau)^\gamma} r(\tau)}{(t - \tau)^\gamma} \xdif \tau \in \xR.
        \end{equation}
        Recall that the set $G_1$ consists of pairs $(t, r(\cdot))$ such that $t \in [0, T]$ and $r(\cdot) \in \C([0, T], \xR)$, and it is endowed with the metric $\dist$ from \eqref{dist}.

        \begin{lmm} \label{lemma_V^1_gamma_mu}
            For every $\gamma \in (0, 1)$ and $\mu > 0$, the following statements hold:
            \begin{itemize}
                \item[$(V.1)$]
                    The functional $V_{\gamma, \mu}$ is continuous.
                \item[$(V.2)$]
                    If $r(\cdot) \in \AC^\gamma([0, T], \xR)$ and $r(0) = 0$, then the function $v(t) = V_{\gamma, \mu}(t, r_t(\cdot))$, $t \in [0, T]$, satisfies the inclusion $v(\cdot) \in \Lip([0, T], \xR)$, and, for a.e. $t \in [0, T]$,
                    \begin{equation*}
                        \dot{v}(t)
                        = (^C D^\gamma r)(t)
                        - \frac{\mu}{\Gamma(1 - \gamma)} r(t)
                        + \frac{\mu^2 \gamma^2}{\Gamma(1 - \gamma)} \int_{0}^{t} \frac{r(\tau)}{(t - \tau)^{\gamma + 1}}
                        \int_{0}^{t - \tau} \xi^{2 \gamma - 1} e^{- \mu \xi^\gamma} \xdif \xi \xdif \tau.
                    \end{equation*}
                    If, in addition, the function $r(\cdot)$ is nonnegative, then
                    \begin{equation*}
                        \dot{v}(t)
                        \leq (^C D^\gamma r)(t)
                        - \frac{\mu}{\Gamma(1 - \gamma)} r(t)
                        + \frac{\mu^2 \Gamma(\gamma + 1)}{2 \Gamma(1 - \gamma)} (I^\gamma r)(t)
                        \text{ for a.e. } t \in [0, T].
                \end{equation*}
            \end{itemize}
        \end{lmm}
        \begin{proof}
            In the proof, we denote $V = V_{\gamma, \mu}$ for brevity.

            1.
            Let us show that, for any $(t, r(\cdot)) \in G_1$, the function $v(\tau) = V(\tau, r_\tau(\cdot))$, $\tau \in [0, t]$, satisfies the estimate
            \begin{equation} \label{proposition_aux_continuity_main}
                |v(\tau^\prime) - v(\tau)|
                \leq \frac{\|r(\cdot)\|_{[0, t]}}{\Gamma(2 - \gamma)} |\tau^\prime - \tau|^{1 - \gamma}
                + \frac{T^{1 - \gamma}}{\Gamma(2 - \gamma)} \varkappa(|\tau^\prime - \tau|),
                \quad \tau, \tau^\prime \in [0, t],
            \end{equation}
            where $\varkappa$ is the modulus of continuity of $r(\cdot)$ on $[0, t]$.

            If $t = 0$, inequality \eqref{proposition_aux_continuity_main} holds automatically.
            So, let $t > 0$.
            Note that
            \begin{equation} \label{change_of_variable}
                v(\tau)
                = \frac{1}{\Gamma(1 - \gamma)} \int_{0}^{\tau} \frac{e^{- \mu \xi^\gamma} r(\tau - \xi)}{\xi^\gamma} \xdif \xi,
                \quad \tau \in [0, t].
            \end{equation}
            Fix $\tau$, $\tau^\prime \in [0, t]$ such that $\tau^\prime > \tau$.
            If $\tau = 0$, then, taking into account that $v(0) = 0$, we obtain
            \begin{equation*}
                |v(\tau^\prime) - v(\tau)|
                = |v(\tau^\prime)|
                \leq \frac{1}{\Gamma(1 - \gamma)} \int_{0}^{\tau^\prime} \frac{e^{- \mu \xi^\gamma} |r(\tau^\prime - \xi)|}{\xi^\gamma} \xdif \xi
                \leq \frac{\|r(\cdot)\|_{[0, t]}}{\Gamma(1 - \gamma)} \int_{0}^{\tau^\prime} \frac{\xdif \xi}{\xi^\gamma}
                = \frac{\|r(\cdot)\|_{[0, t]}}{\Gamma(2 - \gamma)} (\tau^\prime - \tau)^{1 - \gamma},
            \end{equation*}
            and, if $\tau > 0$, we derive
            \begin{align*}
                |v(\tau^\prime) - v(\tau)|
                & \leq \frac{1}{\Gamma(1 - \gamma)} \int_{\tau}^{\tau^\prime} \frac{e^{- \mu \xi^\gamma} |r(\tau^\prime - \xi)|}{\xi^\gamma} \xdif \xi
                + \frac{1}{\Gamma(1 - \gamma)} \int_{0}^{\tau}
                \frac{e^{- \mu \xi^\gamma} |r(\tau^\prime - \xi) - r(\tau - \xi)|}{\xi^\gamma} \xdif \xi \\
                & \leq \frac{\|r(\cdot)\|_{[0, t]}}{\Gamma(1 - \gamma)} \int_{\tau}^{\tau^\prime} \frac{\xdif \xi}{\xi^\gamma}
                + \frac{\varkappa(\tau^\prime - \tau)}{\Gamma(1 - \gamma)} \int_{0}^{\tau} \frac{\xdif \xi}{\xi^\gamma}
                = \frac{\|r(\cdot)\|_{[0, t]}}{\Gamma(2 - \gamma)} \big( (\tau^\prime)^{1 - \gamma} - \tau^{1 - \gamma} \big)
                + \frac{\varkappa(\tau^\prime - \tau)}{\Gamma(2 - \gamma)} \tau^{1 - \gamma} \\
                & \leq \frac{\|r(\cdot)\|_{[0, t]}}{\Gamma(2 - \gamma)} (\tau^\prime - \tau)^{1 - \gamma}
                + \frac{T^{1 - \gamma}}{\Gamma(2 - \gamma)} \varkappa(\tau^\prime - \tau).
            \end{align*}
            Thus, inequality \eqref{proposition_aux_continuity_main} is valid.

            2.
            Now, let us prove statement $(V.1)$.
            Let $(t_0, r^{(0)}(\cdot)) \in G_1$ and $\{(t_k, r^{(k)}(\cdot))\}_{k \in \xN} \subset G_1$ be such that $\dist_k = \dist((t_0, r^{(0)}(\cdot)), (t_k, r^{(k)}(\cdot))) \to 0$ as $k \to \infty$.
            For every $k \in \xN \cup \{0\}$, let $\varkappa_k$ be the modulus of continuity of $r^{(k)}(\cdot)$ on $[0, t_k]$.
            Since the functions $r^{(k)}(\cdot)$, $k \in \xN \cup \{0\}$, are uniformly bounded and equicontinuous (see Sect.~\ref{section_Preliminaries}), there exists $R > 0$ such that $\|r^{(k)}(\cdot)\|_{[0, t_k]} \leq R$, $k \in \xN \cup \{0\}$, and $\varkappa_\ast(\delta) = \sup\{ \varkappa_k(\delta): \, k \in \xN \cup \{0\} \} \to 0$ as $\delta \to + 0$.
            Hence, in order to establish the required convergence $V(t_k, r^{(k)}(\cdot)) \to V(t_0, r^{(0)}(\cdot))$ as $k \to \infty$, it suffices to prove for every $k \in \xN$ the inequality
            \begin{equation} \label{proposition_aux_functional_continuity_proof_first}
                |V(t_0, r^{(0)}(\cdot)) - V(t_k, r^{(k)}(\cdot))|
                \leq \frac{R}{\Gamma(2 - \gamma)} \dist_k^{1 - \gamma}
                + \frac{T^{1 - \gamma}}{\Gamma(2 - \gamma)} (\dist_k + 2 \varkappa_\ast(\dist_k)).
            \end{equation}

            Fix $k \in \xN$.
            Assume that $t_0 \leq t_k$.
            Then, we have
            \begin{equation*}
                |V(t_0, r^{(0)}(\cdot)) - V(t_k, r^{(k)}(\cdot))|
                \leq |V(t_0, r^{(0)}(\cdot)) - V(t_0, r^{(k)}_{t_0}(\cdot))|
                + |V(t_0, r^{(k)}_{t_0}(\cdot)) - V(t_k, r^{(k)}(\cdot))|.
            \end{equation*}
            For the first term, by virtue of \eqref{dist_properties} and \eqref{change_of_variable}, we derive
            \begin{align*}
                |V(t_0, r^{(0)}(\cdot)) - V(t_0, r^{(k)}_{t_0}(\cdot))|
                & \leq \frac{1}{\Gamma(1 - \gamma)} \int_{0}^{t_0} \frac{e^{- \mu \xi^\gamma} |r^{(0)}(t_0 - \xi) - r^{(k)}(t_0 - \xi)|}{\xi^\gamma} \xdif \xi \\
                & \leq \frac{T^{1 - \gamma}}{\Gamma(2 - \gamma)} \max_{\xi \in [0, t_0]}  |r^{(0)}(\xi) - r^{(k)}(\xi)|
                \leq \frac{T^{1 - \gamma}}{\Gamma(2 - \gamma)} (\dist_k + \varkappa_\ast(\dist_k)),
            \end{align*}
            and for the second term, in view of \eqref{dist_properties} and \eqref{proposition_aux_continuity_main}, we get
            \begin{align*}
                |V(t_0, r^{(k)}_{t_0}(\cdot)) - V(t_k, r^{(k)}(\cdot))|
                & \leq \frac{R}{\Gamma(2 - \gamma)} (t_k - t_0)^{1 - \gamma} + \frac{T^{1 - \gamma}}{\Gamma(2 - \gamma)} \varkappa_\ast(t_k - t_0) \\
                & \leq \frac{R}{\Gamma(2 - \gamma)} \dist_k^{1 - \gamma} + \frac{T^{1 - \gamma}}{\Gamma(2 - \gamma)} \varkappa_\ast(\dist_k).
            \end{align*}
            Thus, we obtain inequality \eqref{proposition_aux_functional_continuity_proof_first}.
            In the case $t_k < t_0$, this inequality can be proved in a similar way.

            3.
            Further, let us prove $(V.2)$.
            Fix $r(\cdot) \in \AC^\gamma([0, T], \xR)$ such that $r(0) = 0$ and consider the function $v(t) = V(t, r_t(\cdot))$, $t \in [0, T]$.
            For every $\theta \geq 0$, based on the equality
            \begin{equation} \label{exp}
                e^{- \mu \theta^\gamma}
                = 1 - \mu \gamma \int_{0}^{\theta} \frac{e^{- \mu \xi^\gamma}}{\xi^{1 - \gamma}} \xdif \xi,
            \end{equation}
            which can be verified by direct calculation, we derive
            \begin{equation*}
                e^{- \mu \theta^\gamma}
                = 1 - \mu \gamma \int_{0}^{\theta} \frac{\xdif \xi}{\xi^{1 - \gamma}}
                - \mu \gamma \int_{0}^{\theta} \frac{e^{- \mu \xi^\gamma} - 1}{\xi^{1 - \gamma}} \xdif \xi
                = 1 - \mu \theta^\gamma + \mu \gamma \int_{0}^{\theta} \frac{1 - e^{- \mu \xi^\gamma}}{\xi^{1 - \gamma}} \xdif \xi,
            \end{equation*}
            and, consequently, according to \eqref{I} and \eqref{V^1_gamma_mu}, the function $v(\cdot)$ can be represented as follows:
            \begin{align*}
                v(t)
                & = (I^{1 - \gamma} r)(t)
                - \frac{\mu}{\Gamma(1 - \gamma)} (I^1 r)(t)
                + \frac{\mu \gamma}{\Gamma(1 - \gamma)} \int_{0}^{t} \frac{r(\tau)}{(t - \tau)^\gamma}
                \int_{0}^{t - \tau} \frac{1 - e^{- \mu \xi^\gamma}}{\xi^{1 - \gamma}} \xdif \xi \xdif \tau \\
                & = v_1(t)
                - \frac{\mu}{\Gamma(1 - \gamma)} v_2(t)
                + \frac{\mu \gamma}{\Gamma(1 - \gamma)} v_3(t),
                \quad t \in [0, T].
            \end{align*}
            Since $r(\cdot) \in \AC^\gamma([0, T], \xR)$ and $r(0) = 0$, for $v_1(t) = (I^{1 - \gamma} r)(t)$, $t \in [0, T]$, we conclude $v_1(\cdot) \in \Lip([0, T], \xR)$ and $\dot{v}_1(t) = (^C D^\gamma r)(t)$ for a.e. $t \in [0, T]$ (see Sect.~\ref{section_Preliminaries}).
            For $v_2(t) = (I^1 r)(t)$, $t \in [0, T]$, we have $v_2(\cdot) \in \Lip([0, T], \xR)$, and $\dot{v}_2(t) = r(t)$, $t \in (0, T)$.
            Thus, it remains to investigate the properties of the function
            \begin{equation*}
                v_3(t)
                = \int_{0}^{t} \frac{r(\tau)}{(t - \tau)^\gamma}
                \int_{0}^{t - \tau} \frac{1 - e^{- \mu \xi^\gamma}}{\xi^{1 - \gamma}} \xdif \xi \xdif \tau,
                \quad t \in [0, T].
            \end{equation*}

            4.
            To this end, let us introduce the auxiliary function
            \begin{equation*}
                M(\theta)
                = \frac{1}{\theta^\gamma} \int_{0}^{\theta} \frac{1 - e^{- \mu \xi^\gamma}}{\xi^{1 - \gamma}} \xdif \xi,
                \quad \theta > 0,
            \end{equation*}
            and describe some of its properties.
            It follows from \eqref{exp} that
            \begin{equation*}
                0
                \leq
                1 - e^{- \mu \xi^\gamma}
                \leq \mu \gamma \int_{0}^{\xi} \frac{\xdif \eta}{\eta^{1 - \gamma}}
                = \mu \xi^\gamma,
                \quad \xi \geq 0,
            \end{equation*}
            and, hence,
            \begin{equation} \label{proposition_second_function_proof_estimate_K_1}
                0
                \leq M(\theta)
                \leq \frac{\mu}{\theta^\gamma} \int_{0}^{\theta} \xi^{2 \gamma - 1} \xdif \xi
                = \frac{\mu}{2 \gamma} \theta^{\gamma},
                \quad \theta > 0.
            \end{equation}
            In particular, we obtain $M(\theta) \to 0$ as $\theta \to + 0$.
            Further, by virtue of the integration by parts formula, we derive
            \begin{equation} \label{proposition_second_function_proof_K_dot}
                \dot{M}(\theta)
                = - \frac{\gamma}{\theta^{\gamma + 1}} \int_{0}^{\theta} \frac{1 - e^{- \mu \xi^\gamma}}{\xi^{1 - \gamma}} \xdif \xi
                + \frac{1 - e^{- \mu \theta^\gamma}}{\theta}
                = \frac{\mu \gamma}{\theta^{\gamma + 1}} \int_{0}^{\theta} \xi^{2 \gamma - 1} e^{- \mu \xi^\gamma} \xdif \xi,
                \quad \theta > 0.
            \end{equation}
            Consequently, we have
            \begin{equation} \label{proposition_second_function_proof_estimate_K_2}
                0
                \leq \dot{M}(\theta)
                \leq \frac{\mu \gamma}{\theta^{\gamma + 1}} \int_{0}^{\theta} \xi^{2 \gamma - 1} \xdif \xi
                = \frac{\mu}{2 \theta^{1 - \gamma}},
                \quad \theta > 0,
            \end{equation}
            and, therefore, for any $\theta > 0$ and $\theta^\prime > \theta$,
            \begin{equation} \label{proposition_second_function_proof_estimate_K_3}
                0
                \leq M(\theta^\prime) - M(\theta)
                = \int_{\theta}^{\theta^\prime} \dot{M}(\xi) \xdif \xi
                \leq \int_{\theta}^{\theta^\prime} \frac{\mu}{2 \xi^{1 - \gamma}} \xdif \xi
                \leq \frac{\mu}{2 \theta^{1 - \gamma}} (\theta^\prime - \theta).
            \end{equation}

            5.
            Now, based on the representation
            \begin{equation*}
                v_3(t)
                = \int_{0}^{t} M(t - \tau) r(\tau) \xdif \tau,
                \quad t \in [0, T],
            \end{equation*}
            let us prove first that the function $v_3(\cdot)$ satisfies the Lipschitz condition $|v_3(t^\prime) - v_3(t)| \leq L |t^\prime - t|$, $t$, $t^\prime \in [0, T]$, with the constant
            \begin{equation*}
                L
                = \frac{(\gamma + 2) \mu \|r(\cdot)\|_{[0, T]} T^\gamma}{2 (\gamma + 1) \gamma}.
            \end{equation*}
            Fix $t$, $t^\prime \in [0, T]$ such that $t^\prime > t$.
            If $t = 0$, then, taking into account that $v_3(0) = 0$ and using \eqref{proposition_second_function_proof_estimate_K_1}, we derive
            \begin{align*}
                |v_3(t^\prime) - v_3(t)|
                & = |v_3(t^\prime)|
                \leq \int_{0}^{t^\prime} M(t^\prime - \tau) |r(\tau)| \xdif \tau
                \leq \|r(\cdot)\|_{[0, T]} \int_{0}^{t^\prime} \frac{\mu}{2 \gamma} (t^\prime - \tau)^\gamma \xdif \tau \\
                & = \frac{\mu \|r(\cdot)\|_{[0, T]} (t^\prime)^{\gamma + 1}}{2 (\gamma + 1) \gamma}
                \leq \frac{\mu \|r(\cdot)\|_{[0, T]} T^\gamma}{2 (\gamma + 1) \gamma} t^\prime
                \leq L (t^\prime - t).
            \end{align*}
            Suppose that $t > 0$.
            Then, we have
            \begin{equation} \label{psi_3_Lip_1}
                v_3(t^\prime) - v_3(t)
                = \int_{t}^{t^\prime} M(t^\prime - \tau) r(\tau) \xdif \tau
                + \int_{0}^{t} (M(t^\prime - \tau) - M(t - \tau)) r(\tau) \xdif \tau.
            \end{equation}
            For the first term, according to \eqref{proposition_second_function_proof_estimate_K_1}, we get
            \begin{align}
                \Big| \int_{t}^{t^\prime} M(t^\prime - \tau) r(\tau) \xdif \tau \Big|
                & \leq \int_{t}^{t^\prime} M(t^\prime - \tau) |r(\tau)| \xdif \tau
                \leq \|r(\cdot)\|_{[0, T]} \int_{t}^{t^\prime} \frac{\mu}{2 \gamma} (t^\prime - \tau)^\gamma \xdif \tau \nonumber \\
                & = \frac{\mu \|r(\cdot)\|_{[0, T]} (t^\prime - t)^{\gamma + 1}}{2 (\gamma + 1) \gamma}
                \leq \frac{\mu \|r(\cdot)\|_{[0, T]} T^\gamma}{2 (\gamma + 1) \gamma} (t^\prime - t), \label{proposition_second_function_proof_estimate_I_1}
            \end{align}
            and, for the second term, by virtue of \eqref{proposition_second_function_proof_estimate_K_3}, we conclude
            \begin{align*}
                \Big| \int_{0}^{t} (M(t^\prime - \tau) - M(t - \tau)) r(\tau) \xdif \tau \Big|
                & \leq \int_{0}^{t} (M(t^\prime - \tau) - M(t - \tau)) |r(\tau)| \xdif \tau \nonumber \\
                & \leq \|r(\cdot)\|_{[0, T]} \int_{0}^{t} \frac{\mu}{2 (t - \tau)^{1 - \gamma}} (t^\prime - t) \xdif \tau
                = \frac{\mu \|r(\cdot)\|_{[0, T]} t^\gamma}{2 \gamma} (t^\prime - t) \nonumber \\
                & \leq \frac{\mu \|r(\cdot)\|_{[0, T]} T^\gamma}{2 \gamma} (t^\prime - t).
            \end{align*}
            Thus, we obtain the desired estimate.

            6.
            Since $v_3(\cdot) \in \Lip([0, T], \xR)$, then the derivative $\dot{v}_3(t)$ exists for a.e. $t \in [0, T]$.
            In order to obtain an explicit formula for this derivative, let us calculate the right-hand side derivative $\dot{v}_3^{+}(t)$ of $v_3(\cdot)$ at every $t \in (0, T)$.
            For the first term in \eqref{psi_3_Lip_1}, owing to \eqref{proposition_second_function_proof_estimate_I_1}, we have
            \begin{equation*}
                \Big| \frac{1}{t^\prime - t} \int_{t}^{t^\prime} M(t^\prime - \tau) r(\tau) \xdif \tau \Big|
                \leq \frac{\mu \|r(\cdot)\|_{[0, T]} (t^\prime - t)^\gamma}{2 (\gamma + 1) \gamma},
                \quad t^\prime \in (t, T],
            \end{equation*}
            and, therefore,
            \begin{equation*}
                \lim_{t^\prime \to t + 0} \frac{1}{t^\prime - t} \int_{t}^{t^\prime} M(t^\prime - \tau) r(\tau) \xdif \tau
                = 0.
            \end{equation*}
            Let us consider the second term in \eqref{psi_3_Lip_1}.
            For any $\tau \in [0, t)$, we get
            \begin{equation*}
                \lim_{t^\prime \to t + 0} \frac{(M(t^\prime - \tau) - M(t - \tau)) r(\tau)}{t^\prime - t}
                = \dot{M}(t - \tau) r(\tau),
            \end{equation*}
            and, moreover, due to \eqref{proposition_second_function_proof_estimate_K_3},
            \begin{equation*}
                \frac{|(M(t^\prime - \tau) - M(t - \tau)) r(\tau)|}{t^\prime - t}
                \leq \frac{\mu \|r(\cdot)\|_{[0, T]}}{2 (t - \tau)^{1 - \gamma}},
                \quad t^\prime \in (t, T].
            \end{equation*}
            Then, applying Lebesgue's dominated convergence theorem, we conclude
            \begin{equation*}
                \lim_{t^\prime \to t + 0}
                \frac{1}{t^\prime - t} \int_{0}^{t} \big( M(t^\prime - \tau) - M(t - \tau) \big) r(\tau) \xdif \tau
                = \int_{0}^{t} \dot{M}(t - \tau) r(\tau) \xdif \tau.
            \end{equation*}
            Hence, we derive
            \begin{equation*}
                \dot{v}_3^{+}(t)
                = \lim_{t^\prime \to t + 0} \frac{v_3(t^\prime) - v_3(t)}{t^\prime - t}
                = \int_{0}^{t} \dot{M}(t - \tau) r(\tau) \xdif \tau.
            \end{equation*}
            As a result, in view of \eqref{proposition_second_function_proof_K_dot}, we get
            \begin{equation} \label{psi_3_dot}
                \dot{v}_3(t)
                = \int_{0}^{t} \dot{M}(t - \tau) r(\tau) \xdif \tau
                = \mu \gamma \int_{0}^{t} \frac{r(\tau)}{(t - \tau)^{\gamma + 1}}
                \int_{0}^{t - \tau} \xi^{2 \gamma - 1} e^{- \mu \xi^\gamma} \xdif \xi \xdif \tau
                \text{ for a.e. } t \in [0, T].
            \end{equation}

            7.
            Summarizing the above, we obtain that $v(\cdot) \in \Lip([0, T], \xR)$ and
            \begin{align*}
                \dot{v}(t)
                & = \dot{v}_1(t)
                - \frac{\mu}{\Gamma(1 - \gamma)} \dot{v}_2(t)
                + \frac{\mu \gamma}{\Gamma(1 - \gamma)} \dot{v}_3(t) \\
                & = (^C D^\gamma r)(t)
                - \frac{\mu}{\Gamma(1 - \gamma)} r(t)
                + \frac{\mu^2 \gamma^2}{\Gamma(1 - \gamma)} \int_{0}^{t} \frac{r(\tau)}{(t - \tau)^{\gamma + 1}}
                \int_{0}^{t - \tau} \xi^{2 \gamma - 1} e^{- \mu \xi^\gamma} \xdif \xi \xdif \tau
                \text{ for a.e. } t \in [0, T].
            \end{align*}

            8.
            If the function $r(\cdot)$ is nonnegative, then, for a.e. $t \in [0, T]$, according to \eqref{I}, \eqref{proposition_second_function_proof_estimate_K_2}, and \eqref{psi_3_dot}, we derive
            \begin{equation*}
                \dot{v}_3(t)
                \leq \int_{0}^{t} \frac{\mu r(\tau)}{2 (t - \tau)^{1 - \gamma}} \xdif \tau
                = \frac{\mu \Gamma(\gamma)}{2} (I^\gamma r)(t)
            \end{equation*}
            and, therefore,
            \begin{equation*}
                \dot{v}(t)
                \leq (^C D^\gamma r)(t)
                - \frac{\mu}{\Gamma(1 - \gamma)} r(t)
                + \frac{\mu^2 \Gamma(\gamma + 1)}{2 \Gamma(1 - \gamma)} (I^\gamma r)(t).
            \end{equation*}
            This completes the proof of the lemma.
        \end{proof}

        \subsubsection{Functional $V^\ast_{\beta, \mu}$}

        Let $\beta \in [0, 1 - \alpha)$ and $\mu > 0$ be fixed.
        Note that $\gamma = \alpha + \beta \in (0, 1)$ and take the corresponding functional $V_{\alpha + \beta, \mu}$ from \eqref{V^1_gamma_mu}.
        For every $(t, w(\cdot)) \in G_n$, denote
        \begin{equation*}
            q(\tau \mid t, w(\cdot))
            = \|w(\tau) - w(0)\|^2,
            \quad
            r(\tau \mid t, w(\cdot), \beta)
            = \big(I^\beta q(\cdot \mid t, w(\cdot))\big)(\tau),
            \quad \tau \in [0, t].
        \end{equation*}
        Consider the functional
        \begin{equation} \label{V^2_beta_mu}
            G_n \ni (t, w(\cdot)) \mapsto V^\ast_{\beta, \mu}(t, w(\cdot))
            = V_{\alpha + \beta, \mu} \big(t, r(\cdot \mid t, w(\cdot), \beta)\big) \in \xR.
        \end{equation}
        In accordance with the introduced notations, this functional can be defined explicitly by
        \begin{equation*}
            V^\ast_{\beta, \mu}(t, w(\cdot))
            = \frac{1}{\Gamma(1 - \alpha - \beta) \Gamma(\beta)} \int_{0}^{t}
            \frac{e^{- \mu (t - \tau)^{\alpha + \beta}}}{(t - \tau)^{\alpha + \beta}}
            \int_{0}^{\tau} \frac{\|w(\xi) - w(0)\|^2}{(\tau - \xi)^{1 - \beta}} \xdif \xi \xdif \tau,
            \quad (t, w(\cdot)) \in G_n,
        \end{equation*}
        if $\beta > 0$, and, if $\beta = 0$, by
        \begin{equation} \label{V^2_beta_mu_explicit_beta=0}
            V^\ast_{0, \mu}(t, w(\cdot))
            = \frac{1}{\Gamma(1 - \alpha)} \int_{0}^{t}
            \frac{e^{- \mu (t - \tau)^\alpha} \|w(\tau) - w(0)\|^2}{(t - \tau)^\alpha} \xdif \tau,
            \quad (t, w(\cdot)) \in G_n.
        \end{equation}

        \begin{lmm} \label{lemma_V_beta_mu}
            For every $\beta \in [0, 1 - \alpha)$ and $\mu > 0$, the following statements hold:
            \begin{itemize}
                \item[$(V^\ast.1)$]
                    The functional $V^\ast_{\beta, \mu}$ is continuous.

                \item[$(V^\ast.2)$]
                    For any $(t, w(\cdot)) \in G_n$, the inequality below is valid:
                    \begin{equation*}
                        V^\ast_{\beta, \mu}(t, w(\cdot))
                        \geq e^{- \mu T^{\alpha + \beta}} \big(I^{1 - \alpha} q(\cdot \mid t, w(\cdot))\big)(t).
                    \end{equation*}
                    In particular, the functional $V^\ast_{\beta, \mu}$ is nonnegative.

                \item[$(V^\ast.3)$]
                    If $x(\cdot) \in \AC^\alpha([0, T], \xR^n)$, then the function $v^\ast(t) = V^\ast_{\beta, \mu}(t, x_t(\cdot))$, $t \in [0, T]$, satisfies the inclusion $v^\ast(\cdot) \in \Lip([0, T], \xR)$, and
                    \begin{equation} \label{lemma_V_beta_mu_main}
                        \dot{v}^\ast(t)
                        \leq (^C D^\alpha q)(t)
                        - \frac{\mu}{\Gamma(1 - \alpha - \beta)} (I^\beta q)(t)
                        + \frac{\mu^2 \Gamma(\alpha + \beta + 1)}{2 \Gamma(1 - \alpha - \beta)} (I^{\alpha + 2 \beta} q)(t)
                        \text{ for a.e. } t \in [0, T],
                    \end{equation}
                    where $q(\cdot) = q(\cdot \mid T, x(\cdot))$.
            \end{itemize}
        \end{lmm}
        \begin{proof}
            For brevity, denote $V = V_{\alpha + \beta, \mu}$ and $V^\ast = V^\ast_{\beta, \mu}$.

            1.
            Since the functional $V$ is continuous by $(V.1)$, in order to establish $(V^\ast.1)$, it is sufficient to prove continuity of the mapping
            \begin{equation} \label{auxiliary_mapping}
                G_n \ni (t, w(\cdot)) \mapsto \big(t, r(\cdot \mid t, w(\cdot), \beta) \big) \in G_1.
            \end{equation}
            Take $(t_0, w_0(\cdot)) \in G_n$ and $\{(t_k, w_k(\cdot))\}_{k \in \xN} \subset G_n$ such that $\dist_k = \dist((t_0, w_0(\cdot)), (t_k, w_k(\cdot))) \to 0$ as $k \to \infty$.
            Then, in accordance with Sect.~\ref{section_Preliminaries}, there exists $R > 0$ such that $\|w_k(\cdot)\|_{[0, t_k]} \leq R$ for every $k \in \xN \cup \{0\}$, and $\varkappa_\ast(\delta) = \sup\{ \varkappa_k(\delta): \, k \in \xN \cup \{0\} \} \to 0$ as $\delta \to + 0$, where $\varkappa_k$ is the modulus of continuity  of $w_k(\cdot)$ on $[0, t_k]$.
            Denote $q^{(k)}(\cdot) = q(\cdot \mid t_k, w_k(\cdot))$, $k \in \xN \cup \{0\}$.
            For any $k \in \xN \cup \{0\}$, we have
            \begin{equation*}
                |q^{(k)}(\tau) - q^{(k)}(\tau^\prime)|
                \leq \|w_k(\tau) - w_k(\tau^\prime)\| (\|w_k(\tau)\| + \|w_k(\tau^\prime) + 2 \|w_k(0)\|)
                \leq 4 R \varkappa_\ast(|\tau - \tau^\prime|),
                \quad \tau, \tau^\prime \in [0, t_k].
            \end{equation*}
            Moreover, for every $k \in \xN$, by virtue of \eqref{dist_properties}, we derive
            \begin{equation} \label{q^0_q^k}
                |q^{(0)}(\tau) - q^{(k)}(\tau)|
                \leq 4 R (\|w_0(\tau) - w_k(\tau)\| + \|w_0(0) - w_k(0)\|)
                \leq 8 R (\dist_k + \varkappa_\ast(\dist_k)),
                \quad \tau \in [0, \min\{t_0, t_k\}],
            \end{equation}
            and, hence,
            \begin{equation*}
                \dist\big( (t_0, q^{(0)}(\cdot)), (t_k, q^{(k)}(\cdot)) \big)
                \leq \dist_k + 4 R \varkappa_\ast(\dist_k) + 8 R (\dist_k + \varkappa_\ast(\dist_k)).
            \end{equation*}
            Thus, the sequence $\{(t_k, q^{(k)}(\cdot))\}_{k \in \xN} \subset G_1$ converges to $(t_0, q^{(0)}(\cdot)) \in G_1$.
            Further, consider the functions $r^{(k)}(\cdot) = r(\cdot \mid t_k, w_k(\cdot), \beta)$, $k \in \xN \cup \{0\}$.
            If $\beta = 0$, then $r^{(k)}(\cdot) = q^{(k)}(\cdot)$, $k \in \xN \cup \{0\}$, and, consequently, we get $(t_k, r^{(k)}(\cdot)) \to (t_0, r^{(0)}(\cdot))$ as $k \to \infty$, which proves continuity of mapping \eqref{auxiliary_mapping}.
            Let $\beta > 0$.
            Then, for any $k \in \xN \cup \{0\}$, taking into account the estimate
            \begin{equation*}
                |q^{(k)}(\tau)|
                \leq (\|w_k(\tau)\| + \|w_k(0)\|)^2
                \leq 4 R^2,
                \quad \tau \in [0, t_k],
            \end{equation*}
            and due to \eqref{I_properties}, we obtain
            \begin{equation*}
                |r^{(k)}(\tau) - r^{(k)}(\tau^\prime)|
                \leq \frac{2 |\tau - \tau^\prime|^\beta}{\Gamma(\beta + 1)} \max_{\xi \in [0, t_k]} |q^{(k)}(\xi)|
                \leq \frac{8 R^2}{\Gamma(\beta + 1)} |\tau - \tau^\prime|^\beta,
                \quad \tau, \tau^\prime \in [0, t_k].
            \end{equation*}
            In addition, for any $k \in \xN$, based on \eqref{q^0_q^k}, we derive
            \begin{equation*}
                |r^{(0)}(\tau) - r^{(k)}(\tau)|
                \leq \frac{\tau^\beta}{\Gamma(\beta + 1)} \max_{\xi \in [0, \tau]} |q^{(0)}(\xi) - q^{(k)}(\xi)|
                \leq \frac{8 R T^\beta}{\Gamma(\beta + 1)} (\dist_k + \varkappa_\ast(\dist_k)),
                \quad \tau \in [0, \min\{t_0, t_k\}],
            \end{equation*}
            and, therefore, in view of \eqref{dist_properties}, we have
            \begin{equation*}
                \dist\big( (t_0, r^{(0)}(\cdot)), (t_k, r^{(k)}(\cdot)) \big)
                \leq \dist_k + \frac{8 R^2}{\Gamma(\beta + 1)} \dist_k^\beta + \frac{8 R T^\beta}{\Gamma(\beta + 1)} (\dist_k + \varkappa_\ast(\dist_k)).
            \end{equation*}
            So, the sequence $\{(t_k, r^{(k)}(\cdot))\}_{k \in \xN} \subset G_1$ converges to $(t_0, r^{(0)}(\cdot)) \in G_1$, and, hence, mapping \eqref{auxiliary_mapping} is continuous.

            2.
            Further, for every $(t, w(\cdot)) \in G_n$, since the functions $q(\cdot) = q(\cdot  \mid t, w(\cdot))$ and $r(\cdot \mid t, w(\cdot), \beta) = (I^\beta q)(\cdot)$ are nonnegative, in accordance with \eqref{I} and \eqref{I_semigroup_property}, we obtain
            \begin{align*}
                V^\ast(t, w(\cdot))
                & = \frac{1}{\Gamma(1 - \alpha - \beta)} \int_{0}^{t}
                \frac{e^{- \mu (t - \tau)^{\alpha + \beta}} (I^\beta q)(\tau)}{(t - \tau)^{\alpha + \beta}} \xdif \tau
                \geq \frac{e^{- \mu T^{\alpha + \beta}}}{\Gamma(1 - \alpha - \beta)} \int_{0}^{t}
                \frac{(I^\beta q)(\tau)}{(t - \tau)^{\alpha + \beta}} \xdif \tau \\
                & = e^{- \mu T^{\alpha + \beta}} \big(I^{1 - \alpha - \beta} (I^\beta q)\big)(t)
                = e^{- \mu T^{\alpha + \beta}} (I^{1 - \alpha} q)(t)
                \geq 0,
            \end{align*}
            which proves $(V^\ast.2)$.

            3.
            Let us prove statement $(V^\ast.3)$.
            Fix $x(\cdot) \in \AC^\alpha([0, T], \xR^n)$.
            By virtue of \cite[Corollary~4.2]{Gomoyunov_2018_FCAA}, we have $q(\cdot) = q(\cdot \mid T, x(\cdot)) \in \AC^\alpha([0, T], \xR)$, and, therefore, $r(\cdot) = r(\cdot \mid T, x(\cdot), \beta) = (I^\beta q)(\cdot) \in \AC^{\alpha + \beta}([0, T], \xR)$ in view of the equality $q(0) = 0$.
            Since $r(0) = 0$ and the function $r(\cdot)$ is nonnegative, it follows from $(V.2)$ that the function $v^\ast(t) = V^\ast(t, x_t(\cdot)) = V(t, r_t(\cdot))$, $t \in [0, T]$, satisfies the inclusion $v^\ast(\cdot) \in \Lip([0, T], \xR)$, and
            \begin{equation} \label{lemma_V_beta_mu_p_2}
                \dot{v}^\ast(t)
                \leq (^C D^{\alpha + \beta} r)(t)
                - \frac{\mu}{\Gamma(1 - \alpha - \beta)} r(t)
                + \frac{\mu^2 \Gamma(\alpha + \beta + 1)}{2 \Gamma(1 - \alpha - \beta)} (I^{\alpha + \beta} r)(t)
                \text{ for a.e. } t \in [0, T].
            \end{equation}
            Note that, due to \eqref{I_semigroup_property}, we derive $(I^{\alpha + \beta} r)(t) = (I^{\alpha + 2 \beta} q)(t)$, $t \in [0, T]$, and, moreover, in accordance with \eqref{Caputo},
            \begin{equation*}
                (^C D^{\alpha + \beta} r)(t)
                = \xDrv{(I^{1 - \alpha - \beta} r)(t)}{t}
                = \xDrv{\big( I^{1 - \alpha - \beta} (I^\beta q)\big)(t)}{t}
                = \xDrv{(I^{1 - \alpha} q)(t)}{t}
                = (^C D^\alpha q)(t)
                \text{ for a.e. } t \in [0, T].
            \end{equation*}
            Thus, inequality \eqref{lemma_V_beta_mu_p_2} implies estimate \eqref{lemma_V_beta_mu_main}.
            The lemma is proved.
        \end{proof}

        \subsubsection{Functional $V_\ast$}

        By suitably combining the functionals $V^\ast_{\beta, \mu}$ from \eqref{V^2_beta_mu} for various values of $\beta \in [0, 1 - \alpha)$ and $\mu > 0$, we obtain the following result.
        \begin{lmm} \label{lemma_V}
            For any $\lambda > 0$, there exist a number $\lambda_\ast > 0$ and a functional $V_\ast: G_n \to \xR$ such that:
            \begin{itemize}
                \item[$(V_\ast.1)$]
                    The functional $V_\ast$ is nonnegative and continuous.
                    In addition, if $(t, w(\cdot)) \in G_n$ and $w(\tau) = w(0)$, $\tau \in [0, t]$, then $V_\ast(t, w(\cdot)) = 0$.

                \item[$(V_\ast.2)$]
                    For every function $x(\cdot) \in \AC^\alpha([0, T], \xR^n)$, the function $v_\ast(t) = V_\ast(t, x_t(\cdot))$, $t \in [0, T]$, satisfies the inclusion $v_\ast(\cdot) \in \Lip([0, T], \xR)$, and
                    \begin{equation} \label{lemma_V_main}
                        \dot{v}_\ast(t)
                        \leq e^{- \lambda_\ast t}
                        \big( 2 \langle x(t) - x(0), (^C D^\alpha x)(t) \rangle - \lambda \|x(t) - x(0)\|^2 \big)
                        \text{ for a.e. } t \in [0, T].
                    \end{equation}

                \item[$(V_\ast.3)$]
                    For any compact set $X \subset \C([0, T], \xR^n)$ and any number $\rho > 0$, there exists $\delta > 0$ such that, for every $x(\cdot) \in X$, the inequality $V_\ast(T, x(\cdot)) \leq \delta$ implies the estimate $\|x(\cdot) - x(0)\|_{[0, T]} \leq \rho$.
            \end{itemize}
        \end{lmm}

        \begin{rmrk}
            In the case $\alpha = 1$, statements $(V_\ast.1)$ and $(V_\ast.2)$ are satisfied for the functional
            \begin{equation*}
                V_\ast(t, w(\cdot))
                = e^{- \lambda t} \|w(t) - w(0)\|^2,
                \quad (t, w(\cdot)) \in G_n.
            \end{equation*}
            Thus, Lemma~\ref{lemma_V} provides a suitable analogue of this functional for the case when $\alpha \in (0, 1)$.
        \end{rmrk}

        Before proving Lemma~\ref{lemma_V}, we present an auxiliary proposition.
        \begin{prpstn} \label{proposition_integrals_alpha_beta}
            If $\beta \geq 1 - \alpha$ and a function $\psi(\cdot) \in \xLinfty([0, T], \xR)$ is nonnegative, then
            \begin{equation*}
                (I^\beta \psi)(t)
                \leq \frac{\Gamma(1 - \alpha) T^{\alpha + \beta - 1}}{\Gamma(\beta)} (I^{1 - \alpha} \psi)(t),
                \quad t \in [0, T].
            \end{equation*}
        \end{prpstn}
        \begin{proof}
            According to \eqref{I}, for every $t \in [0, T]$, we have
            \begin{equation*}
                (I^\beta \psi)(t)
                = \frac{1}{\Gamma(\beta)} \int_{0}^{t} \frac{(t - \tau)^{\alpha + \beta - 1} \psi(\tau)}{(t - \tau)^\alpha} \xdif \tau
                \leq \frac{T^{\alpha + \beta - 1}}{\Gamma(\beta)} \int_{0}^{t} \frac{\psi(\tau)}{(t - \tau)^\alpha} \xdif \tau
                \leq \frac{\Gamma(1 - \alpha) T^{\alpha + \beta - 1}}{\Gamma(\beta)} (I^{1 - \alpha} \psi)(t),
            \end{equation*}
            which proves the proposition.
        \end{proof}

        \begin{proof}[Proof of Lemma~\ref{lemma_V}.]
            Fix $\lambda > 0$.
            Choice of the number $\lambda_\ast$ and construction of the functional $V_\ast: G_n \to \xR$ depend on the value of $\alpha$.
            For the reader's convenience, we first consider in detain the cases when $\alpha \in [1 / 2, 1)$ and $\alpha \in [1 / 4, 1 / 2)$, and, after that, we handle the general case when $\alpha \in [2^{- m}, 2^{- (m - 1)})$ for some $m \in \xN$.

            1.
            Assume that $\alpha \in [1 / 2, 1)$.
            Define
            \begin{equation*} \label{case_1_parameters}
                \beta_1
                = 0,
                \quad \mu_1
                = \Gamma(1 - \alpha) \lambda,
                \quad \lambda_\ast
                = \frac{\mu_1^2 \Gamma(\alpha + 1) T^{2 \alpha - 1}}{2 \Gamma(\alpha)} e^{\mu_1 T^\alpha},
            \end{equation*}
            take the corresponding functional $V^\ast_{\beta_1, \mu_1}$ from \eqref{V^2_beta_mu}, and put
            \begin{equation*}
                V_\ast(t, w(\cdot))
                = e^{- \lambda_\ast t} V^\ast_{\beta_1, \mu_1} (t, w(\cdot)),
                \quad (t, w(\cdot)) \in G_n.
            \end{equation*}
            Let us show that the specified $\lambda_\ast$ and $V_\ast$ possess properties $(V_\ast.1)$--$(V_\ast.3)$.

            Since the functional $V^\ast_{\beta_1, \mu_1}$ is nonnegative and continuous by $(V^\ast.1)$ and $(V^\ast.2)$, we obtain that the functional $V_\ast$ is nonnegative and continuous, too.
            Now, let $(t, w(\cdot)) \in G_n$ be such that $w(\tau) = w(0)$, $\tau \in [0, t]$.
            Then, it follows from \eqref{V^2_beta_mu_explicit_beta=0} that $V^\ast_{\beta_1, \mu_1} (t, w(\cdot)) = 0$, and, consequently, $V_\ast(t, w(\cdot)) = 0$.
            Thus, statement $(V_\ast.1)$ is proved.

            Further, fix $x(\cdot) \in \AC^\alpha([0, T], \xR^n)$.
            Introduce the auxiliary function $v_1^\ast(t) = V^\ast_{\beta_1, \mu_1}(t, x_t(\cdot))$, $t \in [0, T]$.
            Due to $(V^\ast.3)$, we derive $v_1^\ast(\cdot) \in \Lip([0, T], \xR)$, and, since $\beta_1 = 0$,
            \begin{equation*}
                \dot{v}^\ast_1(t)
                \leq (^C D^\alpha q)(t)
                - \frac{\mu_1}{\Gamma(1 - \alpha)} q(t)
                + \frac{\mu_1^2 \Gamma(\alpha + 1)}{2 \Gamma(1 - \alpha)} (I^\alpha q)(t)
                \text{ for a.e. } t \in [0, T],
            \end{equation*}
            where $q(\cdot) = q(\cdot \mid T, x(\cdot))$.
            Note that, in the considered case, we have $\alpha \geq 1 - \alpha$.
            Hence, by Proposition~\ref{proposition_integrals_alpha_beta}, taking into account $(V^\ast.2)$, we get
            \begin{equation*}
                (I^\alpha q)(t)
                \leq \frac{\Gamma(1 - \alpha) T^{2 \alpha - 1}}{\Gamma(\alpha)} (I^{1 - \alpha} q)(t)
                \leq \frac{\Gamma(1 - \alpha) T^{2 \alpha - 1}}{\Gamma(\alpha)} e^{\mu_1 T^\alpha} v_1^\ast(t),
                \quad t \in [0, T].
            \end{equation*}
            Therefore, by virtue of the choice of $\mu_1$ and $\lambda_\ast$, we obtain
            \begin{equation*}
                \dot{v}_1^\ast(t)
                \leq (^C D^\alpha q)(t)
                - \lambda q(t)
                + \lambda_\ast v_1^\ast(t)
                \text{ for a.e. } t \in [0, T].
            \end{equation*}
            Thus, for the function $v_\ast(t) = V_\ast(t, x_t(\cdot)) = e^{- \lambda_\ast t} v_1^\ast(t)$, $t \in [0, T]$, we conclude $v_\ast(\cdot) \in \Lip([0, T], \xR)$ and
            \begin{equation*}
                \dot{v}_\ast(t)
                = e^{- \lambda_\ast t} \dot{v}^\ast_1(t) - \lambda_\ast e^{- \lambda_\ast t} v_1^\ast(t)
                \leq e^{- \lambda_\ast t} ((^C D^\alpha q)(t) - \lambda q(t))
                \text{ for a.e. } t \in [0, T].
            \end{equation*}
            From this estimate and the inequality
            \begin{equation} \label{s_derivarive}
                (^C D^\alpha q)(t)
                \leq 2 \langle x(t) - x(0), (^C D^\alpha x)(t) \rangle
                \text{ for a.e. } t \in [0, T],
            \end{equation}
            which is valid by \cite[Corollary~4.2]{Gomoyunov_2018_FCAA}, we derive \eqref{lemma_V_main}.
            Property $(V_\ast.2)$ is established.

            Let us prove $(V_\ast.3)$.
            Arguing by contradiction, suppose that there exist a compact set $X \subset \C([0, T], \xR^n)$ and a number $\rho > 0$ such that, for every $k \in \xN$, one can choose $x^{(k)}(\cdot) \in X$ such that $V_\ast(T, x^{(k)}(\cdot)) \leq 1 / k$ and
            \begin{equation} \label{case_1_jjj}
                \|x^{(k)}(\cdot) - x^{(k)}(0)\|_{[0, T]}
                \geq \rho.
            \end{equation}
            Owing to compactness of $X$, we can assume that the sequence $\{x^{(k)}(\cdot)\}_{k \in \xN}$ converges to a function $x^{(0)}(\cdot) \in X$.
            Denote $q^{(k)}(\cdot) = q(\cdot \mid T, x^{(k)}(\cdot))$, $k \in \xN \cup \{0\}$.
            Then, the sequence $\{q^{(k)}(\cdot)\}_{k \in \xN} \subset \C([0, T], \xR)$ converges to $q^{(0)}(\cdot)$ (see, e.g., the proof of $(V^\ast.1)$), and, therefore, in view of \eqref{I_properties}, we have $(I^{1 - \alpha} q^{(k)})(T) \to (I^{1 - \alpha} q^{(0)})(T)$ as $k \to \infty$.
            On the other hand, for every $k \in \xN$, according to $(V^\ast.2)$, we get
            \begin{equation*}
                1 / k
                \geq V_\ast(T, x^{(k)}(\cdot))
                = e^{- \lambda_\ast T} V^\ast_{\beta_1, \mu_1} (T, x^{(k)}(\cdot))
                \geq e^{- \lambda_\ast T - \mu_1 T^\alpha} (I^{1 - \alpha} q^{(k)})(T)
                \geq 0,
            \end{equation*}
            wherefrom it follows that $(I^{1 - \alpha} q^{(k)})(T) \to 0$ as $k \to \infty$, and, hence, $(I^{1 - \alpha} q^{(0)})(T) = 0$.
            Since the function $q^{(0)}(\cdot)$ is continuous and nonnegative, this equality yields $q^{(0)}(t) = 0$, $t \in [0, T]$.
            But, passing to the limit as $k \to \infty$ in inequality \eqref{case_1_jjj}, we derive $\|q^{(0)}(\cdot)\|_{[0, T]} \geq \rho^2 > 0$, and obtain the contradiction, which proves $(V_\ast.3)$.

            2.
            Consider the next case when $\alpha \in [1 / 4, 1 / 2)$.
            Define
            \begin{equation*}
                \beta_1
                = 0,
                \quad \beta_2
                = \alpha,
                \quad \mu_1
                = 2 \Gamma(1 - \alpha) \lambda,
                \quad \mu_2
                = \frac{\mu_1^2 \Gamma(\alpha + 1) \Gamma(1 - 2 \alpha)}{2 \Gamma(1 - \alpha)}
            \end{equation*}
            and
            \begin{equation*}
                \lambda_\ast
                = \frac{\mu_2^2 \Gamma(2 \alpha + 1) \Gamma(1 - \alpha) T^{4 \alpha - 1}}{2 \Gamma(1 - 2 \alpha) \Gamma(3 \alpha)} e^{\mu_1 T^\alpha}.
            \end{equation*}
            Note that $\beta_2 < 1 - \alpha$ since $\alpha < 1 / 2$.
            Take the functionals $V^\ast_{\beta_1, \mu_1}$ and $V^\ast_{\beta_2, \mu_2}$ from \eqref{V^2_beta_mu} and put
            \begin{equation*}
                V_\ast(t, w(\cdot))
                = e^{- \lambda_\ast t} \big( V^\ast_{\beta_1, \mu_1} (t, w(\cdot)) + V^\ast_{\beta_2, \mu_2} (t, w(\cdot)) \big) / 2,
                \quad (t, w(\cdot)) \in G_n.
            \end{equation*}

            For the specified $\lambda_\ast$ and $V_\ast$, statement $(V_\ast.1)$ is verified in the same way as in the first case.
            The proof of $(V_\ast.3)$ also does not differ essentially from the arguments given above, because, due to $(V^\ast.2)$ and the equality $\beta_1 = 0$, for any $x(\cdot) \in \C([0, T], \xR^n)$, the following estimate holds:
            \begin{equation*}
                V_\ast(T, x(\cdot))
                \geq e^{- \lambda_\ast T} V^\ast_{\beta_1, \mu_1} (T, x(\cdot)) / 2
                \geq e^{- \lambda_\ast T - \mu_1 T^\alpha} (I^{1 - \alpha} q)(T) / 2,
            \end{equation*}
            where $q(\cdot) = q(\cdot \mid T, x(\cdot))$.
            Thus, it remains to prove $(V_\ast.2)$.

            For a given $x(\cdot) \in \AC^\alpha([0, T], \xR^n)$, introduce the functions $v_1^\ast(t) = V^\ast_{\beta_1, \mu_1}(t, x_t(\cdot))$, $v_2^\ast(t) = V^\ast_{\beta_2, \mu_2}(t, x_t(\cdot))$, $t \in [0, T]$.
            According to $(V^\ast.3)$, we obtain $v_1^\ast(\cdot)$, $v_2^\ast(\cdot) \in \Lip([0, T], \xR)$, and, taking into account the choice of $\beta_1$, $\beta_2$, $\mu_1$, and $\mu_2$, we derive
            \begin{align*}
                \dot{v}_1^\ast(t) + \dot{v}_2^\ast(t)
                & \leq 2 (^C D^\alpha q)(t)
                - \frac{\mu_1}{\Gamma(1 - \alpha)} q(t)
                + \frac{\mu_1^2 \Gamma(\alpha + 1)}{2 \Gamma(1 - \alpha)} (I^{\alpha} q)(t) \\
                & + (^C D^\alpha q)(t) - \frac{\mu_2}{\Gamma(1 - 2 \alpha)} (I^\alpha q)(t)
                + \frac{\mu_2^2 \Gamma(2 \alpha + 1)}{2 \Gamma(1 - 2 \alpha)} (I^{3 \alpha} q)(t) \\
                & = 2 (^C D^\alpha q)(t) - 2 \lambda q(t)
                + \frac{\mu_2^2 \Gamma(2 \alpha + 1)}{2 \Gamma(1 - 2 \alpha)} (I^{3 \alpha} q)(t)
                \text{ for a.e. } t \in [0, T],
            \end{align*}
            where $q(\cdot) = q(\cdot \mid T, x(\cdot))$.
            Since, in the second case, we have $3 \alpha \geq 1 - \alpha$, then Proposition~\ref{proposition_integrals_alpha_beta} and $(V^\ast.2)$ yield
            \begin{equation*}
                (I^{3 \alpha} q)(t)
                \leq \frac{\Gamma(1 - \alpha) T^{4 \alpha - 1}}{\Gamma(3 \alpha)} (I^{1 - \alpha} q)(t)
                \leq \frac{\Gamma(1 - \alpha) T^{4 \alpha - 1}}{\Gamma(3 \alpha)} e^{\mu_1 T^\alpha} v_1^\ast(t),
                \quad t \in [0, T].
            \end{equation*}
            Therefore, owing to the choice of $\lambda_\ast$, we conclude
            \begin{equation*}
                \dot{v}_1^\ast(t) + \dot{v}_2^\ast(t)
                \leq 2 (^C D^\alpha q)(t) - 2 \lambda q(t) + \lambda_\ast v_1^\ast(t)
                \text{ for a.e. } t \in [0, T].
            \end{equation*}
            Thus, for the function $v_\ast(t) = V_\ast(t, x_t(\cdot)) = e^{- \lambda_\ast t} (v_1^\ast(t) + v_2^\ast(t)) / 2$, $t \in [0, T]$, we get $v_\ast(\cdot) \in \Lip([0, T], \xR)$ and
            \begin{equation*}
                \dot{v}_\ast(t)
                = e^{- \lambda_\ast t} (\dot{v}_1^\ast(t) + \dot{v}_2^\ast(t)) / 2
                - \lambda_\ast e^{- \lambda_\ast t} (v_1^\ast(t) + v_2^\ast(t)) / 2
                \leq e^{- \lambda_\ast t} ((^C D^\alpha q)(t) - \lambda q(t))
                \text{ for a.e. } t \in [0, T].
            \end{equation*}
            This estimate and \eqref{s_derivarive} imply \eqref{lemma_V_main}.
            The proof of $(V_\ast.3)$ is completed.

            3.
            In the general case, choose $m \in \xN$ such that $\alpha \in [2^{- m}, 2^{- (m - 1)})$.
            The cases $m = 1$ and $m = 2$ were considered above, so we can assume that $m > 2$.
            Put
            \begin{equation*}
                \beta_i
                = (2^{i - 1} - 1) \alpha,
                \quad i \in \overline{1, m}.
            \end{equation*}
            Note that $\beta_1 = 0$, and, due to the choice of $m$,
            \begin{equation*}
                0
                \leq
                \beta_i
                \leq (2^{m - 1} - 1) \alpha
                < (2^{m - 1} - 1) 2^{- (m - 1)}
                = 1 - 2^{- (m - 1)}
                < 1 - \alpha,
                \quad i \in \overline{1, m}.
            \end{equation*}
            Further, define numbers $\mu_i > 0$, $i \in \overline{1, m}$, by the following recurrent relations:
            \begin{equation*}
                \mu_1
                = m \Gamma(1 - \alpha) \lambda,
                \quad \mu_{i + 1}
                = \frac{\mu_i^2 \Gamma(\alpha + \beta_i + 1) \Gamma(1 - \alpha - \beta_{i + 1})}{2 \Gamma(1 - \alpha - \beta_i)},
                \quad i \in \overline{1, m - 1},
            \end{equation*}
            and set
            \begin{equation*}
                \lambda_\ast
                = \frac{\mu_m^2 \Gamma(\alpha + \beta_m + 1) \Gamma(1 - \alpha) T^{2 \alpha + 2 \beta_m - 1}}
                {2 \Gamma(1 - \alpha - \beta_m) \Gamma(\alpha + 2 \beta_m)} e^{\mu_1 T^\alpha}.
            \end{equation*}
            Finally, take the corresponding functionals $V^\ast_{\beta_i, \mu_i}$, $i \in \overline{1, m}$, from \eqref{V^2_beta_mu} and put
            \begin{equation*}
                V_\ast(t, w(\cdot))
                = \frac{e^{- \lambda_\ast t}}{m} \sum_{i = 1}^{m} V^\ast_{\beta_i, \mu_i} (t, w(\cdot)),
                \quad (t, w(\cdot)) \in G_n.
            \end{equation*}

            The proofs of properties $(V_\ast.1)$ and $(V_\ast.3)$ for the specified $\lambda_\ast$ and $V_\ast$ are carried out by the same scheme as in the two particular cases considered above, and, therefore, they are omitted.
            Let us prove $(V_\ast.2)$.

            Let $x(\cdot) \in \AC^\alpha([0, T], \xR^n)$ and $q(\cdot) = q(\cdot \mid T, x(\cdot))$.
            By virtue of $(V^\ast.3)$, for every $i \in \overline{1, m}$, the function $v_i^\ast(t) = V^\ast_{\beta_i, \mu_i}(t, x_t(\cdot))$, $t \in [0, T]$, satisfies the inclusion $v_i^\ast(\cdot) \in \Lip([0, T], \xR)$, and, taking into account the choice of $\beta_1$ and $\mu_1$, we derive
            \begin{align}
                \sum_{i = 1}^{m} \dot{v}_i^\ast(t)
                & \leq m (^C D^\alpha q)(t) - m \lambda q(t) \nonumber \\
                & - \sum_{i = 2}^{m} \frac{\mu_i}{\Gamma(1 - \alpha - \beta_i)} (I^{\beta_i} q)(t)
                + \sum_{i = 1}^{m} \frac{\mu_i^2 \Gamma(\alpha + \beta_i + 1)}{2 \Gamma(1 - \alpha - \beta_i)} (I^{\alpha + 2 \beta_i} q)(t)
                \text{ for a.e. } t \in [0, T]. \label{general_case_sum}
            \end{align}
            Let $t \in [0, T]$ be fixed.
            For any $i \in \overline{1, m - 1}$, due to the choice of $\beta_i$, $\beta_{i + 1}$, and $\mu_{i + 1}$, we have
            \begin{equation*}
                \frac{\mu_i^2 \Gamma(\alpha + \beta_i + 1)}{2 \Gamma(1 - \alpha - \beta_i)} (I^{\alpha + 2 \beta_i} q)(t)
                = \frac{\mu_i^2 \Gamma(\alpha + \beta_i + 1)}{2 \Gamma(1 - \alpha - \beta_i)} (I^{\beta_{i + 1}} q)(t)
                = \frac{\mu_{i + 1}}{\Gamma(1 - \alpha - \beta_{i + 1})} (I^{\beta_{i + 1}} q)(t),
            \end{equation*}
            and, consequently, for the last two terms in \eqref{general_case_sum}, we get
            \begin{equation} \label{general_case_1}
                \sum_{i = 1}^{m} \frac{\mu_i^2 \Gamma(\alpha + \beta_i + 1)}{2 \Gamma(1 - \alpha - \beta_i)} (I^{\alpha + 2 \beta_i} q)(t)
                - \sum_{i = 2}^{m} \frac{\mu_i}{\Gamma(1 - \alpha - \beta_i)} (I^{\beta_i} q)(t)
                = \frac{\mu_m^2 \Gamma(\alpha + \beta_m + 1)}{2 \Gamma(1 - \alpha - \beta_m)} (I^{\alpha + 2 \beta_m} q)(t).
            \end{equation}
            Further, owing to the choice of $\beta_m$ and the inequality $\alpha \geq 2^{- m}$, we obtain
            \begin{equation*}
                \alpha + 2 \beta_m
                = (2^m - 1) \alpha
                \geq (2^m - 1) 2^{- m}
                = 1 - 2^{- m}
                \geq 1 - \alpha.
            \end{equation*}
            Hence, according to Proposition~\ref{proposition_integrals_alpha_beta} and $(V^\ast.2)$, since $\beta_1 = 0$, we conclude
            \begin{equation*}
                (I^{\alpha + 2 \beta_m} q)(t)
                \leq \frac{\Gamma(1 - \alpha) T^{2 \alpha + 2 \beta_m - 1}}{\Gamma(\alpha + 2 \beta_m)} (I^{1 - \alpha} q)(t)
                \leq \frac{\Gamma(1 - \alpha) T^{2 \alpha + 2 \beta_m - 1}}{\Gamma(\alpha + 2 \beta_m)} e^{\mu_1 T^\alpha} v_1^\ast(t),
            \end{equation*}
            and, then, in view of the choice of $\lambda_\ast$, we derive
            \begin{equation} \label{general_case_2}
                \frac{\mu_m^2 \Gamma(\alpha + \beta_m + 1)}{2 \Gamma(1 - \alpha - \beta_m)} (I^{\alpha + 2 \beta_m} q)(t)
                \leq \lambda_\ast v_1^\ast(t).
            \end{equation}
            Thus, it follows from \eqref{general_case_sum}, \eqref{general_case_1}, and \eqref{general_case_2} that
            \begin{equation*}
                \sum_{i = 1}^{m} \dot{v}_i^\ast(t)
                \leq m (^C D^\alpha q)(t)
                - m \lambda q(t)
                + \lambda_\ast v_1^\ast(t)
                \text{ for a.e. } t \in [0, T].
            \end{equation*}
            Therefore, we get that the function
            \begin{equation*}
                v_\ast(t)
                = V_\ast(t, x_t(\cdot))
                = \frac{e^{- \lambda_\ast t}}{m} \sum_{i = 0}^{m} v_i^\ast(t),
                \quad t \in [0, T],
            \end{equation*}
            satisfies the inclusion $v_\ast(\cdot) \in \Lip([0, T], \xR)$, and
            \begin{equation*}
                \dot{v}_\ast(t)
                = \frac{e^{- \lambda_\ast t}}{m} \sum_{i = 0}^{m} \dot{v}_i^\ast(t)
                - \frac{\lambda_\ast e^{- \lambda_\ast t}}{m} \sum_{i = 0}^{m} v_i^\ast(t)
                \leq e^{- \lambda_\ast t} ((^C D^\alpha q)(t) - \lambda q(t))
                \text{ for a.e. } t \in [0, T].
            \end{equation*}
            This estimate and \eqref{s_derivarive} yield \eqref{lemma_V_main}, which completes the proof of $(V_\ast.3)$.
            The lemma is proved.
        \end{proof}

        \subsubsection{Functional $\mathcal{V}_\varepsilon$}
        \label{subsubsection_Veps}

        Below, following the scheme from \cite[Sect.~7.5]{Subbotin_1995} and based on Lemma~\ref{lemma_V}, for every sufficiently small $\varepsilon > 0$, we define a functional $\mathcal{V}_\varepsilon: G_n \to \xR$ with a number of prescribed properties, which are close to those listed in \cite[Sect.~5]{Lukoyanov_2004_PMM_Eng} (in this connection, see assumptions $(H.4)^\prime$ in \cite{Crandall_Lions_Ishii_1987}, $(A.4)$ in \cite[Sect.~9.2]{Subbotin_1995}, and $(F.3)$ in \cite{Lukoyanov_2006_IMM_Eng}).

        Let $R > 0$ be fixed, and let $\lambda_H$ be chosen by $R$ according to assumption $(H.3)$.
        Set $\lambda = 4 \lambda_H$ and take the corresponding number $\lambda_\ast$ and functional $V_\ast$ from Lemma~\ref{lemma_V}.
        Choose $\varepsilon_0 > 0$ such that $\varepsilon_0 \leq 2 e^{- (\lambda_H + \lambda_\ast / 2) T}$.
        For any $\varepsilon \in (0, \varepsilon_0]$, define
        \begin{equation*}
            G_n \ni (t, w(\cdot)) \mapsto \mathcal{V}_\varepsilon(t, w(\cdot))
            = \frac{e^{- \lambda_H t}}{\varepsilon} \sqrt{\varepsilon^4 + V_\ast(t, w(\cdot))}
            \in \xR
        \end{equation*}
        and consider also the auxiliary functionals $p_\varepsilon: G_n \to \xR$ and $s_\varepsilon: G_n \to \xR^n$ given by
        \begin{align*}
            p_\varepsilon(t, w(\cdot))
            & = - \frac{\lambda_H e^{- \lambda_H t}}{\varepsilon} \sqrt{\varepsilon^4 + V_\ast(t, w(\cdot))}
            - \frac{2 \lambda_H e^{- (\lambda_H + \lambda_\ast) t}}{\varepsilon} \frac{\|w(t) - w(0)\|^2}
            {\sqrt{\varepsilon^4 + V_\ast(t, w(\cdot))}}, \\
            s_\varepsilon(t, w(\cdot))
            & = \frac{e^{- (\lambda_H + \lambda_\ast) t}}{\varepsilon} \frac{w(t) - w(0)}{\sqrt{\varepsilon^4 + V_\ast(t, w(\cdot))}},
            \quad (t, w(\cdot)) \in G_n.
        \end{align*}

        \begin{lmm} \label{lemma_nu_varepsilon}
            For every $R > 0$ and $\varepsilon \in (0, \varepsilon_0]$, the following statements hold:
            \begin{itemize}
                \item[$(\mathcal{V}.1)$]
                    The functional $\mathcal{V}_\varepsilon$ is nonnegative and continuous.
                    In addition, if $(t, w(\cdot)) \in G_n$ and $w(\tau) = w(0)$, $\tau \in [0, t]$, then $\mathcal{V}_\varepsilon(t, w(\cdot)) \leq \varepsilon$.

                \item[$(\mathcal{V}.2)$]
                    The functionals $p_\varepsilon$ and $s_\varepsilon$ are continuous.
                    Furthermore, for any $(t, w(\cdot))$, $(t, w^\prime(\cdot)) \in G_n$ such that $\|w(\cdot)\|_{[0, t]} \leq R$, $\|w^\prime(\cdot)\|_{[0, t]} \leq R$, and $w(0) = w^\prime(0)$, the inequality below is valid:
                    \begin{equation} \label{C_2_2_main}
                        p_\varepsilon(t, \Delta w(\cdot))
                        + H\big(t, w^\prime(t), s_\varepsilon(t, \Delta w(\cdot))\big)
                        - H\big(t, w(t), s_\varepsilon(t, \Delta w(\cdot))\big)
                        \leq 0,
                    \end{equation}
                    where $\Delta w(\cdot) = w^\prime(\cdot) - w(\cdot)$.

                \item[$(\mathcal{V}.3)$]
                    For every function $x(\cdot) \in \AC^\alpha([0, T], \xR^n)$, the function $v(t) = \mathcal{V}_\varepsilon(t, x_t(\cdot))$, $t \in [0, T]$, satisfies the inclusion $v(\cdot) \in \Lip([0, T], \xR)$, and
                    \begin{equation} \label{C_2_1_main}
                        \dot{v}(t)
                        \leq p_\varepsilon(t, x_t(\cdot))
                        + \langle s_\varepsilon(t, x_t(\cdot)), (^C D^\alpha x)(t) \rangle
                        \text{ for a.e. } t \in [0, T].
                    \end{equation}

                \item[$(\mathcal{V}.4)$]
                    For any compact set $X \subset \AC^\alpha([0, T], \xR^n)$ and any numbers $K > 0$ and $\kappa > 0$, there exists $\varepsilon_\ast \in (0, \varepsilon_0]$ such that, if $\varepsilon \in (0, \varepsilon_\ast]$ and $x(\cdot)$, $x^\prime(\cdot) \in X$ satisfy the relations $x(0) = x^\prime(0)$ and $\mathcal{V}_\varepsilon(T, x^\prime(\cdot) - x(\cdot)) \leq K$, then $|\sigma(x^\prime(\cdot)) - \sigma(x(\cdot))| \leq \kappa$.
            \end{itemize}
        \end{lmm}

        \begin{rmrk}
            If it were shown that, for some $\varepsilon \in (0, \varepsilon_0]$, the functional $\mathcal{V}_\varepsilon$ is $ci$-smooth of the order $\alpha$ (see Sect.~\ref{subsection_ci-derivatives}), then it would follow from \cite[Lemma~9.2]{Gomoyunov_2020_SIAM} that, for every function $x(\cdot) \in \AC^\alpha([0, T], \xR^n)$, the function $v(t) = \mathcal{V}_\varepsilon(t, x_t(\cdot))$, $t \in [0, T]$, would satisfy the equality
            \begin{equation*}
                \dot{v}(t)
                = \partial^\alpha_t \mathcal{V}_\varepsilon(t, x_t(\cdot))
                + \langle \nabla^\alpha \mathcal{V}_\varepsilon(t, x_t(\cdot)), (^C D^\alpha x)(t) \rangle
                \text{ for a.e. } t \in [0, T].
            \end{equation*}
            Comparing this equality with estimate \eqref{C_2_1_main}, we see that, in some sense, the functionals $p_\varepsilon$ and $s_\varepsilon$ correspond to the derivatives $\partial^\alpha_t \mathcal{V}_\varepsilon$ and $\nabla^\alpha \mathcal{V}_\varepsilon$, respectively.
            Hence, in this case, statement $(\mathcal{V}.2)$ could be considered as the requirement for the functionals $\partial^\alpha_t \mathcal{V}_\varepsilon$ and $\nabla^\alpha \mathcal{V}_\varepsilon$, which is consistent with the item $(d)$ in \cite[Sect.~5]{Lukoyanov_2004_PMM_Eng}.
            However, for the results of the paper to be valid, $ci$-smoothness of the order $\alpha$ of the functional $\mathcal{V}_\varepsilon$ is not necessary, and we only need to establish the properties given in Lemma~\ref{lemma_nu_varepsilon}.
        \end{rmrk}

        \begin{proof}[Proof of Lemma~\ref{lemma_nu_varepsilon}.]
            1.
            Taking into account that the functional $V_\ast$ is nonnegative and continuous by $(V_\ast.1)$, we obtain that the functional $\mathcal{V}_\varepsilon$ is nonnegative and continuous, too.
            Now, let $(t, w(\cdot)) \in G_n$ be such that $w(\tau) = w(0)$, $\tau \in [0, t]$.
            Then, in view of $(V_\ast.1)$, we derive $\mathcal{V}_\varepsilon(t, w(\cdot)) = e^{- \lambda_H t} \varepsilon \leq \varepsilon$.
            Thus, property $(\mathcal{V}.1)$ is established.

            2.
            Note that continuity of the functionals $V_\ast$ and $G_n \ni (t, w(\cdot)) \mapsto w(t) - w(0) \in \xR^n$ imply also continuity of the functionals $p_\varepsilon$ and $s_\varepsilon$.
            Further, fix $(t, w(\cdot))$, $(t, w^\prime(\cdot)) \in G_n$ such that $\|w(\cdot)\|_{[0, t]} \leq R$, $\|w^\prime(\cdot)\|_{[0, t]} \leq R$, and $w(0) = w^\prime(0)$ and denote $\Delta w(\cdot) = w^\prime(\cdot) - w(\cdot)$.
            By the choice of $\lambda_H$, we have
            \begin{align} \label{C_2_2_main_p1}
                & p_\varepsilon(t, \Delta w(\cdot))
                + H\big(t, w^\prime(t), s_\varepsilon(t, \Delta w(\cdot))\big)
                - H\big(t, w(t), s_\varepsilon(t, \Delta w(\cdot))\big) \nonumber \\
                & \quad \leq p_\varepsilon(t, \Delta w(\cdot)) + \lambda_H \big(1 + \|s_\varepsilon(t, \Delta w(\cdot))\|\big) \|\Delta w(t)\| \nonumber \\
                & \quad = - \frac{\lambda_H e^{- \lambda_H t}}{\varepsilon} \sqrt{\varepsilon^4 + V_\ast(t, \Delta w(\cdot))}
                \Big( 1 - \varepsilon e^{\lambda_H t} \frac{\|\Delta w(t)\|}{\sqrt{\varepsilon^4 + V_\ast(t, \Delta w(\cdot))}}
                + e^{- \lambda_\ast t} \frac{\|\Delta w(t)\|^2}{\varepsilon^4 + V_\ast(t, \Delta w(\cdot))} \Big).
            \end{align}
            Due to the choice of $\varepsilon_0$, we get $\varepsilon e^{\lambda_H t} \leq \varepsilon_0 e^{\lambda_H T} \leq 2 e^{- \lambda_\ast t / 2}$, and, therefore,
            \begin{equation} \label{C_2_2_main_p2}
                1 - \varepsilon e^{\lambda_H t} \frac{\|\Delta w(t)\|}{\sqrt{\varepsilon^4 + V_\ast(t, \Delta w(\cdot))}}
                + e^{- \lambda_\ast t} \frac{\|\Delta w(t)\|^2}{\varepsilon^4 + V_\ast(t, \Delta w(\cdot))}
                \geq \Big( 1 - e^{- \lambda_\ast t / 2} \frac{\|\Delta w(t)\|}{\sqrt{\varepsilon^4 + V_\ast(t, \Delta w(\cdot))}} \Big)^2
                \geq 0.
            \end{equation}
            From \eqref{C_2_2_main_p1} and \eqref{C_2_2_main_p2}, we derive \eqref{C_2_2_main}.

            3.
            Now, let $x(\cdot) \in \AC^\alpha([0, T], \xR^n)$ and $v(t) = \mathcal{V}_\varepsilon(t, x_t(\cdot))$, $t \in [0, T]$.
            Then, it follows from $(V_\ast.2)$ that $v(\cdot) \in \Lip([0, T], \xR)$ and, by virtue of the choice of $\lambda$,
            \begin{align*}
                \dot{v}(t)
                & = - \frac{\lambda_H e^{- \lambda_H t}}{\varepsilon} \sqrt{\varepsilon^4 + V_\ast(t, x_t(\cdot))}
                + \frac{e^{- \lambda_H t}}{2 \varepsilon \sqrt{\varepsilon^4 + V_\ast(t, x_t(\cdot))}} \xDrv{V_\ast(t, x_t(\cdot))}{t} \\
                & \leq - \frac{\lambda_H e^{- \lambda_H t}}{\varepsilon} \sqrt{\varepsilon^4 + V_\ast(t, x_t(\cdot))}
                + \frac{e^{- (\lambda_H + \lambda_\ast) t}}{\varepsilon \sqrt{\varepsilon^4 + V_\ast(t, x_t(\cdot))}}
                \langle x(t) - x(0), (^C D^\alpha x)(t) \rangle \\
                & - \frac{2 \lambda_H e^{- (\lambda_H + \lambda_\ast) t}}{\varepsilon \sqrt{\varepsilon^4 + V_\ast(t, x_t(\cdot))}}
                \|x(t) - x(0)\|^2
                \text{ for a.e. } t \in [0, T].
            \end{align*}
            This estimate, in accordance with the definitions of the functionals $p_\varepsilon$ and $s_\varepsilon$, yields \eqref{C_2_1_main}.

            4.
            Finally, let us prove $(\mathcal{V}.4)$.
            Let $X \subset \AC^\alpha([0, T], \xR^n)$ be a compact set, and let $K > 0$ and $\kappa > 0$.
            Taking into account that the functional $\sigma$ is continuous by assumption $(\sigma)$, choose $\rho > 0$ such that, for any $x(\cdot)$, $x^\prime(\cdot) \in X$, the inequality $\|x^\prime(\cdot) - x(\cdot)\|_{[0, T]} \leq \rho$ implies the estimate $|\sigma(x^\prime(\cdot)) - \sigma(x(\cdot))| \leq \kappa$.
            Consider the set
            \begin{equation*}
                \Delta X
                = \big\{ \Delta x(\cdot) = x^\prime(\cdot) - x(\cdot):
                \, x(\cdot), x^\prime(\cdot) \in X \big\}
                \subset \AC^\alpha([0, T], \xR^n).
            \end{equation*}
            Since $\Delta X$ is compact, based on $(V_\ast.3)$, take $\delta > 0$ such that, for every $\Delta x(\cdot) \in \Delta X$, it follows from the inequality $V_\ast(T, \Delta x(\cdot)) \leq \delta$ that $\|\Delta x(\cdot) - \Delta x(0)\|_{[0, T]} \leq \rho$.
            Now, choose $\varepsilon_\ast \in (0, \varepsilon_0]$ from the condition $K^2 e^{2 \lambda_H T} \varepsilon_\ast^2 \leq \delta$.
            Let us show that statement $(\mathcal{V}.4)$ is valid for the specified $\varepsilon_\ast$.

            Let $\varepsilon \in (0, \varepsilon_\ast]$ and $x(\cdot)$, $x^\prime(\cdot) \in X$ be fixed such that the function $\Delta x(\cdot) = x^\prime(\cdot) - x(\cdot) \in \Delta X$ satisfies the relations $\Delta x(0) = 0$ and $\mathcal{V}_{\varepsilon}(T, \Delta x(\cdot)) \leq K$.
            Then, we derive
            \begin{equation*}
                V_\ast(T, \Delta x(\cdot))
                \leq \varepsilon^4 + V_\ast(T, \Delta x(\cdot))
                = \varepsilon^2 e^{2 \lambda_H T} \mathcal{V}^2_\varepsilon(T, \Delta x(\cdot))
                \leq \delta,
            \end{equation*}
            and, therefore, we have $\|\Delta x(\cdot)\|_{[0, T]} \leq \rho$, wherefrom we obtain $|\sigma(x^\prime(\cdot)) - \sigma(x(\cdot))| \leq \kappa$.
            The lemma is proved.
        \end{proof}

    \subsection{Proof of Theorem~\ref{theorem_comparison_priciple}}

        1.
        Fix $(t_0, w_0(\cdot)) \in G_n^\alpha$.
        If $t_0 = T$, then inequality \eqref{theorem_comparison_priciple_main} for $(t_0, w_0(\cdot))$ holds due to boundary conditions \eqref{upper_solution_boundary} for $\varphi_+$ and \eqref{lower_solution_boundary} for $\varphi_-$.
        So, we further assume that $t_0 < T$.
        Put
        \begin{equation*}
            X_\ast^\alpha(t_0, w_0(\cdot))
            = \big\{ x(\cdot) \in X^\alpha(t_0, w_0(\cdot)):
            \, \|(^C D^\alpha x)(t)\| \leq c_H (1 + \|x(t)\|)
            \text{ for a.e. } t \in [t_0, T] \big\},
        \end{equation*}
        where the set $X^\alpha(t_0, w_0(\cdot))$ is given by \eqref{X^alpha}, and $c_H$ is the constant from assumption $(H.2)$.
        Owing to Proposition~\ref{proposition_DI_1}, the set $X_\ast^\alpha(t_0, w_0(\cdot))$ is compact.
        In particular, there exists $R > 0$ such that $\|x(\cdot)\|_{[0, T]} \leq R$ for any $x(\cdot) \in X^\alpha(t_0, w_0(\cdot))$.
        By this number $R$, define the number $\varepsilon_0 > 0$ and the functionals $\mathcal{V}_\varepsilon$, $p_\varepsilon$, and $s_\varepsilon$ for every $\varepsilon \in (0, \varepsilon_0]$ according to Sect.~\ref{subsubsection_Veps}.

        2.
        Let $\varepsilon \in (0, \varepsilon_0]$.
        For any $t \in [0, T]$ and $\bar{w}(\cdot) = (w(\cdot), w^\prime(\cdot)) \in \AC^\alpha([0, t], \xR^n \times \xR^n)$, consider the set
        \begin{align}
            & \mathcal{F}_\varepsilon(t, w(\cdot), w^\prime(\cdot)) \nonumber \\
            & \quad = \Big\{ (f, f^\prime, h) \in \xR^n \times \xR^n \times \xR: \,
            \|f\| \leq c_H (1 + \|w(t)\|),
            \, \|f^\prime\| \leq c_H (1 + \|w^\prime(t)\|), \nonumber \\
            & \qquad \qquad \qquad \quad \big|h - \langle s_\varepsilon (t, \Delta w(\cdot)), f - f^\prime \rangle
            - H\big(t, w^\prime(t), s_\varepsilon (t, \Delta w(\cdot))\big)
            + H\big(t, w(t), s_\varepsilon (t, \Delta w(\cdot))\big) \big|
            \leq \varepsilon \Big\}, \label{F_varepsilon}
        \end{align}
        where we denote $\Delta w(\cdot) = w^\prime(\cdot) - w(\cdot)$.
        Thus, in accordance with notation \eqref{G^alpha_n}, we obtain the set-valued functional $G_{2 n}^\alpha \ni (t, \bar{w}(\cdot) = (w(\cdot), w^\prime(\cdot))) \mapsto \mathcal{F}_\varepsilon(t, w(\cdot), w^\prime(\cdot)) \subset \xR^{2 n} \times \xR$.
        Note that $\mathcal{F}_\varepsilon$ possesses properties $(\mathcal{F}.1)$--$(\mathcal{F}.3)$.
        Indeed, $(\mathcal{F}.1)$ and $(\mathcal{F}.3)$ can be verified directly, and $(\mathcal{F}.2)$ follows from continuity of the function $H$ (see $(H.1)$) and the functional $s_\varepsilon$ (see $(\mathcal{V}.2)$).
        Further, take $z_0 = \varphi_+(t_0, w_0(\cdot)) - \varphi_-(t_0, w_0(\cdot))$ and consider the Cauchy problem for the functional differential inclusion
        \begin{equation} \label{DI_proof}
            \big( (^C D^\alpha x)(t), (^C D^\alpha x^\prime)(t), \dot{z}(t) \big)
            \in \mathcal{F}_\varepsilon(t, x_t(\cdot), x^\prime_t(\cdot)),
            \quad (x(t), x^\prime(t), z(t)) \in \xR^n \times \xR^n \times \xR, \quad t \in [t_0, T],
        \end{equation}
        and the initial condition
        \begin{equation} \label{DI_proof_initial_condition}
            x(t)
            = x^\prime(t)
            = w_0(t),
            \quad z(t)
            = z_0,
            \quad t \in [0, t_0].
        \end{equation}
        By Proposition~\ref{proposition_FDI}, the set $\mathcal{W}_\varepsilon$ of solutions $(x(\cdot), x^\prime(\cdot), z(\cdot))$ of problem \eqref{DI_proof}, \eqref{DI_proof_initial_condition} is nonempty and compact in $\C([0, T], \xR^n \times \xR^n \times \xR)$.

        3.
        Let us show that there are functions $(x^{(\varepsilon)}(\cdot), x^{\prime (\varepsilon)}(\cdot), z^{(\varepsilon)}(\cdot)) \in \mathcal{W}_\varepsilon$ satisfying the inequality
        \begin{equation} \label{x^varepsilon}
            z^{(\varepsilon)}(T)
            \geq \varphi_+(T, x^{(\varepsilon)}(\cdot)) - \varphi_-(T, x^{\prime (\varepsilon)}(\cdot)).
        \end{equation}
        For every $t \in [t_0, T]$, consider the set
        \begin{equation} \label{M_varepsilon}
            \mathcal{M}_\varepsilon(t)
            = \big\{ (x(\cdot), x^\prime(\cdot), z(\cdot)) \in \mathcal{W}_\varepsilon:
            \, z(t)
            \geq \varphi_+ (t, x_t(\cdot)) - \varphi_- (t, x^\prime_t(\cdot)) \big\}.
        \end{equation}
        Note that $\mathcal{M}_\varepsilon(t_0) \neq \emptyset$ by virtue of initial condition \eqref{DI_proof_initial_condition} and the choice of $z_0$.
        Put
        \begin{equation} \label{t_varepsilon}
            t_\varepsilon
            = \max\big\{t \in [t_0, T]:
            \, \mathcal{M}_\varepsilon(t) \neq \emptyset \big\}.
        \end{equation}
        The maximum is achieved here owing to compactness of $\mathcal{W}_\varepsilon$, lower semicontinuity of $\varphi_+$, and upper semicontinuity of $\varphi_-$.
        So, in order to complete the proof, it is sufficient to verify that $t_\varepsilon = T$.

        Arguing by contradiction, assume that $t_\varepsilon < T$.
        Take
        \begin{equation} \label{hat_x_x^prime_z}
            (\hat{x}(\cdot), \hat{x}^\prime(\cdot), \hat{z}(\cdot))
            \in \mathcal{M}_\varepsilon(t_\varepsilon)
        \end{equation}
        and denote $\hat{s} = s_\varepsilon (t_\varepsilon, \hat{x}^\prime_{t_\varepsilon}(\cdot) - \hat{x}_{t_\varepsilon}(\cdot))$.
        Since $\varphi_+$ and $\varphi_-$ possess respectively properties $(\varphi_+^\ast)$ and $(\varphi_-^\ast)$, choose characteristics $(x^+(\cdot), z^+(\cdot)) \in CH(t_\varepsilon, \hat{x}_{t_\varepsilon}(\cdot), 0, \hat{s})$ and $(x^-(\cdot), z^-(\cdot)) \in CH(t_\varepsilon, \hat{x}^\prime_{t_\varepsilon}(\cdot), 0, \hat{s})$ such that
        \begin{equation} \label{x^+_x^-}
            \varphi_+(t, x^+_t(\cdot)) - z^+(t)
            \leq \varphi_+(t_\varepsilon, \hat{x}_{t_\varepsilon}(\cdot)),
            \quad \varphi_-(t, x^-_t(\cdot)) - z^-(t)
            \geq \varphi_-(t_\varepsilon, \hat{x}^\prime_{t_\varepsilon}(\cdot)),
            \quad t \in [t_\varepsilon, T].
        \end{equation}
        In accordance with \eqref{E} and \eqref{F_varepsilon}, continuity of $H$ and $s_\varepsilon$ implies that there exists $\delta \in (0, T - t_\varepsilon)$ such that $((^C D^\alpha x^+)(t), (^C D^\alpha x^-)(t), \dot{z}^+(t) - \dot{z}^-(t)) \in \mathcal{F}_\varepsilon(t, x^+_t(\cdot), x^-_t(\cdot))$ for a.e. $t \in [t_\varepsilon, t_\varepsilon + \delta]$.
        Then, in view of \eqref{M_varepsilon}, \eqref{hat_x_x^prime_z}, and \eqref{x^+_x^-}, for the function $\bar{z}(t) = \hat{z}(t_\varepsilon) + z^+(t) - z^-(t)$, $t \in [t_\varepsilon, t_\varepsilon + \delta]$, we get
        \begin{align*}
            \bar{z}(t_\varepsilon + \delta)
            & \geq \varphi_+(t_\varepsilon, \hat{x}_{t_\varepsilon}(\cdot)) - \varphi_-(t_\varepsilon, \hat{x}^\prime_{t_\varepsilon}(\cdot))
            + z^+(t_\varepsilon + \delta) - z^-(t_\varepsilon + \delta) \\
            & \geq \varphi_+(t_\varepsilon + \delta, x^+_{t_\varepsilon + \delta}(\cdot))
            - \varphi_-(t_\varepsilon + \delta, x^-_{t_\varepsilon + \delta}(\cdot)).
        \end{align*}
        Further, let functions $\tilde{x}^\pm(\cdot) \in X^\alpha(t_\varepsilon + \delta, x^\pm_{t_\varepsilon + \delta}(\cdot))$ be such that $(^C D^\alpha \tilde{x}^\pm)(t) = 0$ for a.e. $t \in [t_\varepsilon + \delta, T]$ (in this connection, see, e.g., \cite[Lemma~3]{Gomoyunov_2020_DE}).
        Denote $\tilde{s}(t) = s_\varepsilon(t, \tilde{x}^-_t(\cdot) - \tilde{x}^+_t(\cdot))$, $t \in [t_\varepsilon + \delta, T]$, and consider the function $\tilde{z}: [0, T] \to \xR$ defined by $\tilde{z}(t) = \hat{z}(t)$ for $t \in [0, t_\varepsilon]$, $\tilde{z}(t) = \bar{z}(t)$ for $t \in (t_\varepsilon, t_\varepsilon + \delta]$, and
        \begin{equation*}
            \tilde{z}(t)
            = \bar{z}(t_\varepsilon + \delta) + \int_{t_\varepsilon + \delta}^{t}
            \big( H(\tau, \tilde{x}^-(\tau), \tilde{s}(\tau)) - H(\tau, \tilde{x}^+(\tau), \tilde{s}(\tau)) \big) \xdif \tau,
            \quad t \in (t_\varepsilon + \delta, T].
        \end{equation*}
        Hence, by construction, we obtain $(\tilde{x}^+(\cdot), \tilde{x}^-(\cdot), \tilde{z}(\cdot)) \in \mathcal{M}_\varepsilon(t_\varepsilon + \delta)$, which contradicts definition \eqref{t_varepsilon} of $t_\varepsilon$.

        4.
        Note that, for every $\varepsilon \in (0, \varepsilon_0]$, it follows from the inclusion $(x^{(\varepsilon)}(\cdot), x^{\prime (\varepsilon)}(\cdot), z^{(\varepsilon)}(\cdot)) \in \mathcal{W}_\varepsilon$ that $x^{(\varepsilon)}(\cdot)$, $x^{\prime (\varepsilon)}(\cdot) \in X^\alpha_\ast(t_0, w_0(\cdot))$, $z^{(\varepsilon)}(\cdot) \in \Lip([0, T], \xR)$, $z(t_0) = \varphi_+(t_0, w_0(\cdot)) - \varphi_-(t_0, w_0(\cdot))$, and
        \begin{align*}
            \dot{z}^{(\varepsilon)}(t)
            & \leq \langle s_\varepsilon (t, \Delta x^{(\varepsilon)}_t(\cdot)),
            (^C D^\alpha x^{(\varepsilon)})(t) - (^C D^\alpha x^{\prime (\varepsilon)})(t) \rangle \\
            & + H\big(t, x^{\prime (\varepsilon)}(t), s_\varepsilon (t, \Delta x^{(\varepsilon)}_t(\cdot)\big)
            - H\big(t, x^{(\varepsilon)}(t), s_\varepsilon (t, \Delta x^{(\varepsilon)}_t(\cdot))\big)
            + \varepsilon
            \text{ for a.e. } t \in [t_0, T],
        \end{align*}
        where $\Delta x^{(\varepsilon)}(\cdot) = x^{\prime (\varepsilon)}(\cdot) - x^{(\varepsilon)}(\cdot)$.
        Hence, for the function $v(t) = \mathcal{V}_\varepsilon(t, \Delta x^{(\varepsilon)}_t(\cdot)) + z^{(\varepsilon)}(t) - \varepsilon (t - t_0)$, $t \in [t_0, T]$, due to $(\mathcal{V}.1)$, we have
        \begin{equation*}
            v(t_0)
            = \mathcal{V}_\varepsilon(t_0, w_0(\cdot) - w_0(\cdot)) + z^{(\varepsilon)}(t_0)
            \leq \varepsilon + \varphi_+(t_0, w_0(\cdot)) - \varphi_-(t_0, w_0(\cdot)),
        \end{equation*}
        and, in accordance with $(\mathcal{V}.3)$, we obtain $v(\cdot) \in \Lip([0, T], \xR)$ and
        \begin{align*}
            \dot{v}(t)
            & = \xDrv{\mathcal{V}_\varepsilon(t, \Delta x^{(\varepsilon)}_t(\cdot))}{t} + \dot{z}^{(\varepsilon)}(t) - \varepsilon \\
            & \leq p_\varepsilon(t, \Delta x^{(\varepsilon)}_t(\cdot))
            + H\big(t, x^{\prime (\varepsilon)}(t), s_\varepsilon (t, \Delta x^{(\varepsilon)}_t(\cdot))\big)
            - H\big(t, x^{(\varepsilon)}(t), s_\varepsilon (t, \Delta x^{(\varepsilon)}_t(\cdot))\big)
            \text{ for a.e. } t \in [t_0, T],
        \end{align*}
        wherefrom, by virtue of $(\mathcal{V}.2)$ and the choice of $R$, we derive $\dot{v}(t) \leq 0$ for a.e. $t \in [t_0, T]$.
        Thus, we conclude
        \begin{equation*}
            \mathcal{V}_\varepsilon(T, \Delta x^{(\varepsilon)}(\cdot)) + z^{(\varepsilon)}(T) - \varepsilon (T - t_0)
            = v(T)
            \leq v(t_0)
            \leq \varepsilon + \varphi_+(t_0, w_0(\cdot)) - \varphi_-(t_0, w_0(\cdot)).
        \end{equation*}
        Since \eqref{x^varepsilon} and boundary conditions \eqref{upper_solution_boundary} for $\varphi_+$ and \eqref{lower_solution_boundary} for $\varphi_-$ imply that $z^{(\varepsilon)}(T) \geq \sigma(x^{(\varepsilon)}(\cdot)) - \sigma(x^{\prime (\varepsilon)}(\cdot))$, we finally get the estimate
        \begin{equation} \label{main_estimate}
            \mathcal{V}_\varepsilon(T, \Delta x^{(\varepsilon)}(\cdot)) + \sigma(x^{(\varepsilon)}(\cdot)) - \sigma(x^{\prime (\varepsilon)}(\cdot))
            \leq \varepsilon (1 + T - t_0) + \varphi_+(t_0, w_0(\cdot)) - \varphi_-(t_0, w_0(\cdot)),
        \end{equation}
        which is valid for every $\varepsilon \in (0, \varepsilon_0]$.

        In view of compactness of $X^\alpha_\ast(t_0, w_0(\cdot))$, take $K > 0$ such that $|\sigma(x^{(\varepsilon)}(\cdot)) - \sigma(x^{\prime (\varepsilon)}(\cdot))| \leq K$ for every $\varepsilon \in (0, \varepsilon_0]$.
        Then, due to \eqref{main_estimate}, we have
        \begin{equation*}
            \mathcal{V}_\varepsilon(T, \Delta x^{(\varepsilon)}(\cdot))
            \leq K + \varepsilon_0 (1 + T - t_0) + \varphi_+(t_0, w_0(\cdot)) - \varphi_-(t_0, w_0(\cdot)),
            \quad \varepsilon \in (0, \varepsilon_0],
        \end{equation*}
        and, therefore, applying $(\mathcal{V}.4)$, we obtain $\sigma(x^{(\varepsilon)}(\cdot)) - \sigma(x^{\prime (\varepsilon)}(\cdot)) \to 0$ as $\varepsilon \to + 0$.
        Further, since the functionals $\mathcal{V}_\varepsilon$, $\varepsilon \in (0, \varepsilon_0]$, are nonnegative (see $(\mathcal{V}.1)$), it also follows from \eqref{main_estimate} that
        \begin{equation*}
            \sigma(x^{(\varepsilon)}(\cdot)) - \sigma(x^{\prime (\varepsilon)}(\cdot))
            \leq \varepsilon (1 + T - t_0) + \varphi_+(t_0, w_0(\cdot)) - \varphi_-(t_0, w_0(\cdot)),
            \quad \varepsilon \in (0, \varepsilon_0].
        \end{equation*}
        Passing to the limit as $\varepsilon \to + 0$ in this estimate, we derive inequality \eqref{theorem_comparison_priciple_main} for $(t_0, w_0(\cdot))$.
        The theorem is proved.

\section{Existence and Uniqueness}
\label{section_Existence_Uniqueness}

    The main result of the paper is the following theorem, which is valid under assumptions $(H.1)$--$(H.3)$ and $(\sigma)$ from Sect.~\ref{subsection_HJE}.
    \begin{thrm}
        There exists a unique minimax solution of problem \eqref{HJ}, \eqref{HJ_boundary_condition}.
    \end{thrm}
    \begin{proof}
        Taking into account Theorem~\ref{theorem_comparison_priciple}, in order to prove the statement, it is sufficient to show that there exist an upper solution $\varphi_+^\circ: G_n^\alpha \to \xR$ and a lower solution $\varphi_-^\circ: G_n^\alpha \to \xR$ of problem \eqref{HJ}, \eqref{HJ_boundary_condition} such that $\varphi_+^\circ(t, w(\cdot)) \leq \varphi_-^\circ(t, w(\cdot))$ for every $(t, w(\cdot)) \in G_n^\alpha$.
        Construction of such functionals $\varphi_+^\circ$ and $\varphi_-^\circ$ repeats essentially the arguments given in \cite[Sect.~7]{Lukoyanov_2003_1} and follows the scheme from \cite[Theorem~8.2]{Subbotin_1995} (see also \cite[Sect.~5]{Bayraktar_Keller_2018} and \cite[Theorem~1]{Plaksin_2019_DE_Eng}).
        The basis of this construction are the properties of the set-valued function $E$ from \eqref{E} and the sets of characteristics, which are provided by Propositions~\ref{proposition_DI_1}, \ref{proposition_DI_2}, and~\ref{proposition_semigroup_property}.
        For the reader's convenience, we briefly outline the main steps of the proof below.

        Let $\Phi_+$ be the set of functionals $\varphi: G_n^\alpha \to \xR$ that satisfy boundary condition \eqref{upper_solution_boundary} and possess property $(\varphi_+)$.
        Respectively, by $\Phi_-$, we denote the set of functionals $\varphi: G_n^\alpha \to \xR$ such that \eqref{lower_solution_boundary} and $(\varphi_-)$ are valid.

        1.
        For a given $s \in \xR^n$, consider the functionals $\psi_+^{(s)}: G_n^\alpha \to \xR$ and $\psi_-^{(s)}: G_n^\alpha \to \xR$ defined by
        \begin{equation*}
            \psi_+^{(s)}(t, w(\cdot))
            = \max_{(x(\cdot), z(\cdot)) \in CH(t, w(\cdot), 0, s)} \big( \sigma(x(\cdot)) - z(T) \big),
            \quad \psi_-^{(s)}(t, w(\cdot))
            = \min_{(x(\cdot), z(\cdot)) \in CH(t, w(\cdot), 0, s)} \big( \sigma(x(\cdot)) - z(T) \big),
        \end{equation*}
        where $(t, w(\cdot)) \in G_n^\alpha$.
        The functionals $\psi_+^{(s)}$ and $\psi_-^{(s)}$ are respectively upper and lower semicontinuous, and
        \begin{equation*}
            \psi_+^{(s)}(T, w(\cdot))
            = \psi_-^{(s)}(T, w(\cdot))
            = \sigma(w(\cdot)),
            \quad w(\cdot) \in \AC^\alpha([0, T], \xR^n).
        \end{equation*}
        Further, the inclusion $\psi_+^{(s)} \in \Phi_+$ holds, and, in particular, the set $\Phi_+$ is not empty.
        In addition, for every functional $\varphi \in \Phi_+$, the inequality below is valid:
        \begin{equation*}
            \varphi(t, w(\cdot))
            \geq \psi_-^{(s)}(t, w(\cdot)),
            \quad (t, w(\cdot)) \in G_n^\alpha.
        \end{equation*}

        2.
        Put
        \begin{equation*}
            \varphi^\circ(t, w(\cdot))
            = \inf\big\{ \varphi(t, w(\cdot)):
            \, \varphi \in \Phi_+\big\},
            \quad (t, w(\cdot)) \in G_n^\alpha.
        \end{equation*}
        For any $s \in \xR^n$, we have
        \begin{equation*}
            \psi_-^{(s)}(t, w(\cdot))
            \leq \varphi^\circ(t, w(\cdot))
            \leq \psi_+^{(s)}(t, w(\cdot)),
            \quad (t, w(\cdot)) \in G_n^\alpha,
        \end{equation*}
        and, hence, $\varphi^\circ(T, w(\cdot)) = \sigma(w(\cdot))$, $w(\cdot) \in \AC^\alpha([0, T], \xR^n)$.
        Moreover, the functional $\varphi^\circ: G_n^\alpha \to \xR$ possesses property $(\varphi_+)$, and, consequently, we obtain $\varphi^\circ \in \Phi_+$.

        3.
        For every $\vartheta \in [0, T]$ and $s \in \xR^n$, the functional $\varphi^{(\vartheta, s)}: G_n^\alpha \to \xR$ given by
        \begin{equation*}
            \varphi^{(\vartheta, s)}(t, w(\cdot))
            = \begin{cases}
                \displaystyle
                \sup_{(x(\cdot), z(\cdot)) \in CH(t, w(\cdot), 0, s)} \big( \varphi^\circ(\vartheta, x_\vartheta(\cdot)) - z(\vartheta) \big),
                & \mbox{if } t \in [0, \vartheta), \\
                \varphi^\circ(t, w(\cdot)),
                & \mbox{if } t \in [\vartheta, T],
              \end{cases}
              \quad (t, w(\cdot)) \in G_n^\alpha,
        \end{equation*}
        satisfies the inclusion $\varphi^{(\vartheta, s)} \in \Phi_+$.
        Based on this fact, we derive that $\varphi^\circ \in \Phi_-$.

        4.
        Finally, we define the required functionals $\varphi_+^\circ$ and $\varphi_-^\circ$ as respectively the lower and upper closures of the functional $\varphi^\circ$:
        \begin{equation*}
            \varphi_+^\circ(t, w(\cdot))
            = \lim_{\delta \to + 0} \inf_{(t^\prime, w^\prime(\cdot)) \in O_\delta(t, w(\cdot))} \varphi^\circ(t^\prime, w^\prime(\cdot)),
            \quad \varphi_-^\circ(t, w(\cdot))
            = \lim_{\delta \to + 0} \sup_{(t^\prime, w^\prime(\cdot)) \in O_\delta(t, w(\cdot))} \varphi^\circ(t^\prime, w^\prime(\cdot)),
        \end{equation*}
        where
        \begin{equation*}
            O_\delta(t, w(\cdot))
            = \big\{ (t^\prime, w^\prime(\cdot)) \in G_n^\alpha:
            \, \dist\big( (t, w(\cdot)), (t^\prime, w^\prime(\cdot)) \big) \leq \delta \big\},
            \quad (t, w(\cdot)) \in G_n^\alpha.
        \end{equation*}
        Then, $\varphi_+^\circ$ is an upper solution of problem \eqref{HJ}, \eqref{HJ_boundary_condition}, and $\varphi_-^\circ$ is a lower solution of this problem.
        Moreover, by construction, we obtain $\varphi_+^\circ(t, w(\cdot)) \leq \varphi_-^\circ(t, w(\cdot))$, $(t, w(\cdot)) \in G_n^\alpha$.
        This completes the proof of the theorem.
    \end{proof}

\section{Conclusion}
\label{section_Conclusion}

    In the paper, a Cauchy problem for a Hamilton--Jacobi equation with $ci$-derivatives of an order $\alpha \in (0, 1)$ has been considered.
    A notion of a generalized in the minimax sense solution of this problem has been proposed.
    It has been proved that a minimax solution exists, is unique, and is consistent with the classical solution of the problem.
    A special attention has been given to construction of a suitable Lyapunov--Krasovskii functional needed for the proof of a comparison principle.

    Possible directions for further research in this area include but are not limited to the following:
    \begin{itemize}
        \item[(i)]
            establish a relation between the value functional in an optimal control problems for a dynamical system described by differential equations with the Caputo fractional derivatives and the minimax solution of the corresponding Hamilton--Jacobi--Bellman equation;
            obtain such results for differential games;

        \item[(ii)]
            find an infinitesimal criteria for the minimax solution in terms of suitable directional derivatives
            (see, e.g., \cite[Sect.~6.3]{Subbotin_1995} and also \cite{Lukoyanov_2003_1,Lukoyanov_2006_IMM_Eng,Lukoyanov_Plaksin_2019_MIAN_Eng});

        \item[(iii)]
            develop the theory of generalized in the viscosity sense (see, e.g., \cite{Crandall_Lions_1983} and also \cite{Lukoyanov_2007_IMM_Eng}) solutions of the Cauchy problem considered in the paper.
    \end{itemize}

\bibliographystyle{plain}
\bibliography{references}

\end{document}